\documentclass[article, 11pt]{article}
\usepackage{graphicx}

\usepackage[usenames, dvipsnames]{color}

\usepackage{bbm}

\usepackage{color}
\usepackage{latexsym}
\usepackage{amsmath, amssymb, amsfonts}
\usepackage[toc]{appendix}
\numberwithin{equation}{section}

\RequirePackage[utf8]{inputenc}
\newtheorem{theorem}{Theorem}[section]

\newtheorem{definition}[theorem]{Definition}

\newtheorem{proposition}[theorem]{Proposition}

\newtheorem{assumption}[theorem]{Assumption}

\newtheorem{construction}[theorem]{Construction}

\newtheorem{remark}[theorem]{Remark}
\newtheorem{lemma}[theorem]{Lemma}

\newcommand{\noi}{\noindent}

\newcommand{\E}{\mathbb{E}}
\newcommand{\R}{\mathbb{R}}

\newcommand{\N}{\mathbb{N}}

\newcommand{\la}{\lambda}

\newcommand{\eps}{\varepsilon}

\newcommand{\ph}{\varphi}

\newcommand{\al}{\alpha}

\newcommand{\om}{\omega}

\newcommand{\Gam}{\mathnormal{\Gamma}}

\newcommand{\Om}{\mathnormal{\Omega}}

\newcommand{\bbS}{{\mathbb S}}

\newcommand{\PP}{{\mathbb P}}

\newcommand{\calC}{{\cal C}}
\newcommand{\calD}{{\cal D}}

\newcommand{\calF}{{\cal F}}
\newcommand{\calG}{{\cal G}}
\newcommand{\calH}{{\cal H}}

\newcommand{\calM}{{\cal M}}

\newcommand{\calP}{{\cal P}}

\newcommand{\calPTL}{{\calP_{T,L}}}

\newcommand{\skp}{\vspace{\baselineskip}}

\newcommand\iy{\infty}

\newcommand{\clc}{{\cal C}}
\newcommand{\clb}{{\cal B}}

\newcommand{\clr}{{\cal R}}
\newcommand{\clo}{{\cal O}}
\newcommand{\clm}{{\cal M}}
\newcommand{\cld}{{\cal D}}
\newcommand{\cls}{{\cal S}}
\newcommand{\clp}{{\cal P}}
\newcommand{\clf}{{\cal F}}
\newcommand{\clg}{{\cal G}}
\newcommand{\cla}{{\cal A}}
\newcommand{\one}{\mathbbm{1}}

\title{
A numerical scheme for a mean field game in some
queueing systems based on Markov chain approximation method
\thanks{This is the final version of the paper. To appear in {\it SIAM Journal on Control and Optimization.}}}
\author{Erhan Bayraktar\thanks{Research partially supported by the National Science Foundation (DMS-1613170) and the Susan M.~Smith Professorship.}
\and
Amarjit Budhiraja\thanks{Research partially supported by the National Science Foundation (DMS-1305120), the Army Research Office (W911NF-14-1-0331) and DARPA (W911NF-15-2-0122).}
\and
Asaf Cohen}

\date{\today}

\begin{document}

\maketitle

\begin{abstract}
\noi We use the Markov chain approximation method to construct approximations for the solution of the mean field game (MFG) with reflecting barriers studied in \cite{BBC2017}. The  MFG is formulated in terms of a controlled reflected diffusion with a cost function that depends on the reflection terms in addition to the standard variables: state, control, and the mean field term.  This MFG arises from the asymptotic analysis of an $N$-player 
game for single server queues with strategic servers. By showing that our scheme is an almost contraction, we establish the convergence of this numerical scheme over a small time interval.

\skp

\noi{\bf AMS Classification:} 65M12, 60K25, 91A13, 60K35, 93E20,
65M12, 60F17.
\skp

\noi{\bf Keywords:} Numerical scheme, mean field games, Nash equilibrium, rate control, reflected diffusions, heavy traffic limits, queuing systems, Markov chain approximation method.
\end{abstract}

\section{Introduction}\label{sec1}

The theory of mean field games (MFGs) was initiated a decade ago in the seminal work of Lasry and Lions \cite{Lasry2006,Lasry2006b,Lasry2007}, and Huang, Malham{\'e}, and Caines \cite{Huang2006,Huang2007}. For theoretical study and applications of this theory see \cite{Cardaliaguet2013,Gueant2013,Gomes2013,Carmona2013,Lacker2015,Carmona2015,Lacker2015general,Fischer2016} and the references therein. MFGs are control problems that approximate many player games with weak interaction between the players that is given in terms of the empirical distribution of the players' states. 
In these control problems the empirical distribution that governs the interaction is replaced by a deterministic flow of measures. A solution of the MFG is a probability measure on the path space of the single player state 
that is the distribution of the state process under the optimal control for the control problem associated with the 
flow of measures given by the (time-)marginal distributions of this probability measure.

\skp\noi
A standard (probablistic) method to prove the existence of a MFG solution is by solving a fixed point theorem on the space of probability measures on certain path spaces. 
A probability measure on the path space is fixed and a stochastic control problem is formulated in terms of the flow of time marginals of this probability measure.
Then a `best reply' to the probability measure  is found by solving this  control problem. The distribution of the state process under the best reply is another measure on the path space. 
Solution of the MFG is the fixed point of this map that takes a probability measure to its implied `best reply' distribution.
There are other ways to describe a MFG solution, for example the seminal papers of Lasry and Lions \cite{Lasry2006,Lasry2006b,Lasry2007} represent a MFG solution through two coupled nonlinear partial differential equations; one is an equation of Hamilton-Jacobi-Bellman (HJB) type while the second takes the form of a Kolmogorov forward equation, and  recent works of Carmona, Delarue, and Lacker \cite{Carmona2013, Carmona2015}, using probabilistic methods,  characterize the MFG solution as a solution to certain forward backward stochastic differential equations. 
In general closed form solutions for MFGs are not available and thus one needs numerical approximations. In our work we study one such procedure that uses the Markov chain approximation method (\cite{Kushner1992}) and establish convergence of the scheme over a small time interval.

\skp\noi
In recent years there have been several works on numerical schemes for MFGs, most of which are based on the PDE system of \cite{Lasry2007}. 
Achdou and Capuzzo-Dolcetta \cite{Achdou2010} were the first to suggest a finite difference method for approximating the PDE system relying on monotone approximations of the Hamiltonians and a weak formulation of the forward equation. Together with Camilli, the same authors proved in \cite{Achdou2013} the convergence of the scheme. In \cite{Achdou2016}, Achdou and Porretta showed that the solutions of a certain discrete system converges to a weak solution of the PDE system. In \cite{Lachapelle2010}, Lachapelle, Salomon, and Turinici provided an iterative scheme using a discrete Markov decision problem. Taking advantage of the structure of the problem (in particular, the problem is linear-quadratic in the control), they used the monotonic algorithm method introduced in \cite{Salomon2011} and iteratively constructed a value function, control, and a measure by using finite differences based on the forward-backward system. Gu\'eant studies numerical schemes when the Hamiltonians are quadratic, see \cite{Gueant2012a,Gueant2012}. Semi-Lagrangian schemes were studied by Carlini and Silva in \cite{Carlini2014,Carlini2015}. In a recent paper, Chassagneux, Crisan, and Delarue \cite{Chassagneux2017} used the master equation and by making smoothness assumptions on this infinite dimensional PDE, they proposed an algorithm based on Picard iterations and the continuation method.  
The master equation is a parabolic partial differential equation with a terminal condition. Its variables are time, state, and measure and its solution approximates the value function of the MFG, see e.g., \cite{Carmona2014, CDLL2015, DLR2017, BC2017, CP2017}. 

\skp\noi
Our method in contrast to above methods is purely probabilistic. We do not make smoothness assumptions as in \cite{Chassagneux2017}. We use an iterative Markov chain approximation method (see \cite{Kushner1992}) to construct numerical solutions of the MFG. Specifically, we discretize time and space and for a fixed measure on the path space we define a Markov decision problem that is suggested by the MFG. In the first step of the iteration, the law of the solution of the MDP is computed. Then we take this law as the starting point to formulate the MDP for the second iteration and repeat the process. Unfortunately, it is not clear that the map defined by such iterations is in general a contraction. We instead show that the map is an \emph{almost contraction} over a small time interval with length independent of the discretization parameter. By an almost contraction we roughly mean that the map is a contraction
up to  an additional term  that vanishes as the discretization parameter approaches $0$. 
The proof of this almost contraction property relies on the construction of a coupling between certain controlled reflected Markov chains (see proof of Proposition \ref{prop_h_cont}) which we believe is of independent interest. 
Using the above almost contraction property, tightness of relevant processes, and weak convergence arguments, we show the convergence of the laws obtained from the iteration scheme to the solution of the MFG over a small time interval.
Proving the convergence of a Markov chain based approximation method of the form considered in this work over an arbitrary time interval is for now a challenging open problem. 


\skp\noi
The paper is organized as follows. In Section \ref{sec2} we present the MFG and summarize the results from \cite{BBC2017}. In Section \ref{sec3} we provide the numerical scheme and present our main convergence result (Theorem \ref{thm:thm52}). Section \ref{sec4} provides proofs of some auxiliary results from Section \ref{sec3}. Finally Section \ref{sec:numerics} gives a numerical example.


\subsection{Preliminaries} \label{sec11}
We use the following notation.
For every $t\in(0,\iy)$ and
 $f:[0,\iy)\to\R^d$, let $\|f\|_t\doteq\sup_{[0,t]}\|f\|$. In case that $d=1$, we often use $|f|_t$. 
For any two metric spaces $\cls_1, \cls_2$ denote by $\clc(\cls_1:\cls_2)$ the space of continuous functions mapping $\cls_1$ to $\cls_2$. When $\cls_2 = \R$, we use the notation $\clc(\cls_1)$. For a Polish space $\cls$, the space $\clc([0,T]: \cls)$ will be equipped with the
 uniform topology. We will denote by $\cld([0,T]: \cls)$ the space of functions mapping $[0,T]$ to $\cls$ that are right continuous and have left limits (RCLL) defined on $[0,T]$. This space is equipped with the usual Skorohod topology.
Denote by $\calP(\cls)$ the space of probability measures on $\cls$. We endow $\calP(\cls)$ with the topology of weak convergence of measures. Convergence in distribution of $S$ valued random variable $X_n$ to $X$ will be denoted as
$X_n \Rightarrow X$. For $T, L \in (0,\infty)$, the space $\clp(\clc([0,T]:[0,L]))$ will be denoted as $\clp_{T,L}$.
The Wasserstein distance of order $1$ on $\clp(\cls)$, where $\cls$ is a compact metric space,
 is defined as
\begin{align}\notag
W_1(\eta',\eta)=\inf\left\{\left[\int_{\cls}d(x,y)d\pi(x,y)\right] : \pi\in\calP(\cls\times \cls)\;\;\text{ with marginals $\eta'$ and $\eta$}\right\},
\end{align}
 where $\eta,\eta'\in\calP(\cls)$.
For $\phi \in \calC^{1,2}([0,T]\times [0,L])$, $D_t\phi, D\phi, D^2\phi$ will denote the time derivative and the first two space derivatives of $\phi$, respectively.
For $x \in \cls$, $\delta_x \in \clp(\cls)$ denotes the
Dirac measure at $x$.

Throughout the paper we will make extensive use of the Skorohod map, which for the particular setting of interest here is recalled below.  Fix $T,L>0$.
\begin{definition}\label{def_Skorohod}
Given $\psi\in\cld([0,T]:\R)$ such that $\psi(0)\in [0,L]$, we say the triplet of functions $(\ph,\zeta_1,\zeta_2)\in\calD([0,T]: \R^3)$ solve the Skorohod problem for $\psi$ if the following properties are satisfied:

\noi
(i) For every $t\in[0,T], \;\ph(t)=\psi(t)+\zeta_1(t)-\zeta_2(t)\in [0,L]$.

\noi
(ii) $\zeta_i$ are nonnegative and nondecreasing, $\zeta_1(0)=\zeta_2(0)=0$, and
\begin{align}\notag
\int_{[0,T]}1_{(0,L]}(\ph(s))d\zeta_1(s)=\int_{[0,T]}1_{[0,L)}(\ph(s))d\zeta_2(s)=0.
\end{align}
We denote by $\Gamma(\psi)=(\Gamma_1,\Gamma_2,\Gamma_3)(\psi)\doteq (\ph,\zeta_1,\zeta_2)$ and refer to $\Gamma$ as the Skorohod map.
\end{definition}
It is known that there is a unique solution to the Skorohod problem for every $\psi \in\calD([0,T]: \R)$ and so the Skorohod map in Definition \ref{def_Skorohod} is well defined.
The Skorohod map has the following Lipschitz property (see \cite{Kruk2007}).
\begin{lemma}\label{lem_Skorohod}
There exists $c_S \in (0,\infty)$ such that
for all $\om,\tilde\om\in\calD([0, T]:\R)$ with $\om(0), \tilde \om(0) \in [0,L]$,
\begin{equation}\notag
\sum_{i=1}^3\|\Gam_i(\om) - \Gam_i(\tilde\om)\|_{T}
\le c_S\|\om-\tilde\om\|_T.
\end{equation}
\end{lemma}

\section{The MFG and related results}\label{sec2} 
We now provide a precise description of the MFG that was studied in \cite{BBC2017} and state some relevant results from there. 
\subsection{Description of the MFG}\label{sec2a}
Fix $L,T>0$. Here $T$ denotes the terminal time of our finite time horizon and $[0,L]$ will be the state space of the controlled process $X$. Also, let $U$ be a compact subset of $\R$ representing the control space. Let $(\Omega,\calF,\{\calF_t\},\PP)$ be a filtered probability space that supports a one dimensional standard $\calF_t$-Brownian motion $B$. We will refer to the collection $(\Omega,\calF,\{\calF_t\},\PP, B)$
as a system and denote it by $\Xi$.
Given $(x,t,\nu) \in [0,L]\times [0,T]\times \clp_{T,L}$, we denote by $\cla(\Xi, t,x,\nu)$ the collection of
all pairs $(\alpha, Z)$ where $\alpha = \{\alpha(s)\}_{0\le s \le T-t}$ is a $U$-valued
$\clf_s$-progressively measurable process, $Z = \{Z(s)\}_{0\le s \le T-t}$ is a $[0,L]\times \R_+\times \R_+$ valued $\clf_s$-adapted continuous process such that, $Z=(X,Y,R)$ and
\begin{align}\label{13}
Z(s) = (X,Y,R)(s)=\Gamma\left(x+\int_0^\cdot\bar b(u)du+\sigma B(\cdot)\right)(s),\quad s\in[0,T-t],
\end{align}
where 
\begin{align}\notag
\bar b(u)\doteq b(t+u,\nu(t+u),X(u), \al(u)),\; u \in [0, T-t],
\end{align}
$b:[0,T]\times\clp([0,L])\times[0,L]\times U$, $\nu(s)$ is the marginal of $\nu$ at time instant $s$ and $\sigma$ is a (strictly) positive constant.
Given $\nu \in \clp_{T,L}$, $t \in [0,T]$, $x \in [0,L]$, and a system
$\Xi$ as above, let $(\alpha,Z) \in \cla(\Xi, t,x,\nu)$. The cost functon is given by,
 \begin{align}\notag
J_\nu(t,x,\al,Z)&\doteq \E\Big[\int_0^{T-t} f(s+t,\nu(s+t),X(s),\al(s) )ds + g(\nu(T),X(T-t))\\
&\quad+\int_0^{T-t}y(s+t,\nu(s+t))dY(s)+\int_0^{T-t}r(s+t,\nu(s+t))dR(s)\Big],\label{eq:eq234}
\end{align}
and the value function is
\begin{align}\label{15}
V_\nu(t,x)=\inf_{\Xi}\inf_{(\al,Z)\in \cla(\Xi,t,x,\nu)}J_\nu(t,x,\al,Z).
\end{align}
Conditions on $f,g,y,r$ will be specified below. We now introduce the notion of a solution to the MFG associated with
\eqref{13}--\eqref{15}.
\begin{definition}
	\label{def:mfg729}
A {\bf solution to the MFG}, 	associated with
	\eqref{13}--\eqref{15}, with initial condition $x \in [0,L]$ is defined to be a $\nu \in \clp_{T,L}$ such that
there exist a system $\Xi$ and an $(\alpha,Z) \in \cla(\Xi,0,x,\nu)$ such that $Z=(X,Y,R)$
satisfies $\PP \circ X^{-1} = \nu$ and
\begin{equation}\label{eq:eq716}
	V_\nu(0,x) = J_{\nu}(0, x, \alpha,Z).
\end{equation}
If there exists a unique such $\nu$, we refer to $V_\nu(0,x)$ as the {\bf value of the MFG} with initial condition
$x$.
\end{definition}

\subsection{Background results}
  
The following conditions were used in \cite{BBC2017} in order to characterize  the value function $V_\nu$ and the optimal control.
\begin{assumption}\label{assumption1}$\,$
\begin{enumerate}
\item[(a)] There exists $c_L \in (0,\infty)$ such that for every $(t,\eta,x,\al),(t',\eta',x',\al')\in[0,T]\times\calP([0,L])\times[0,L]\times U$,
\begin{align}
&|f(t,\eta,x,\al)-f(t',\eta',x',\al')| + |g(\eta,x) - g(\eta',x')|+|b(t,\eta,x,\al)-b(t',\eta',x',\al')|\nonumber\\
&\quad +|y(t,\eta)-y(t',\eta')|+|r(t,\eta)-r(t',\eta')|\nonumber\\
&\quad\quad\le c_L (|t-t'|+W_1(\eta,\eta')+|x-x'|+|\al-\al'|).\label{19b}
\end{align}
\item[(b)]
For every $(t,\eta, x,p) \in [0,T]\times \clp([0,L])\times [0,L]\times \R$, there is a unique
$\hat \al(t,\eta,x,p) \in U$ such that
\begin{align}\label{19a}
\hat \al(t,\eta,x,p)=\underset{u\in U}{\arg\min}\;h(t,\eta,x,u,p) ,
\end{align}
where 
\begin{align} \label{19}
h(t,\eta,x,u,p)=f(t,\eta,x,u)+b(t,\eta,x,u)p.
\end{align}

\end{enumerate}
\end{assumption}
As argued in \cite{BBC2017}, Berge's maximum theorem (see \cite[Theorem 17.31]{Aliprantis2006}) together with part (b) of the above assumption
implies that 
$\hat \al$ is continuous.
Also note that \eqref{19b} implies that $b,f,g,y,r$ are bounded functions, in particular,
\begin{equation}\sup_{(t,\eta,x,u)\in [0,T]\times \clp([0,L])\times [0,L]\times U}|b(t,\eta,x,u)| \doteq c_B < \infty.
	\label{eq:eq429}
\end{equation}
For further discussion about the assumption, see \cite{BBC2017}.

The Hamilton-Jacobi-Bellman equation for the value function $V_{\nu}(t,x)$ is given as follows.
\begin{align}\label{HJB1}
-D_t \phi-H(t,\nu(t),x,D \phi)-\frac{1}{2}\sigma^2 D^2\phi=0,\qquad (t,x)\in[0,T]\times[0,L],
\end{align}
with the boundary conditions
\begin{align}\label{HJB2}
\phi(T,x)=g(\nu(T),x),\; D\phi(t,0)=-y(t,\nu(t)),\text{ and } D\phi(t,L)=r(t,\nu(t)),\; t\in[0,T],
\end{align}
where $H$ is the Hamiltonian given as
$$H(t,\eta,x,p)=\inf_{u\in U} h (t,\eta,x,u,p).$$
The following class of H\"{o}lder continuous $\nu \in \clp_{T,L}$ plays a key role in the analysis.
$$\calM_0:=\{\nu \in \clp_{T,L}:\sup_{0\le s<t\le T}\frac{W_1(\nu(t),\nu(s))}{(t-s)^{1/2}}<\iy
\}.$$
\begin{proposition}[Proposition 3.1 in \cite{BBC2017}]\label{lem2}
Fix  $\nu\in \calM_0$ and suppose that Assumption \ref{assumption1} holds. Then $V_\nu$ is continuously differentiable w.r.t.~$t$ and twice continuously differentiable with respect to (w.r.t.) $x$. It is the unique solution of the Hamilton-Jacobi-Bellman 
equation \eqref{HJB1} with 
the boundary conditions \eqref{HJB2}.
Furthermore, with $\hat \alpha$ as introduced in Assumption \ref{assumption1},
%
the map $(s,x')\mapsto \hat \al(s,\nu(s),x',D V_\nu(s,x'))$ is continuous and the
feedback control $\hat\gamma(u, x')\doteq \hat \al(u+t,\nu(u+t),x', D V_\nu(u+t,x'))$ is an optimal feedback control for \eqref{15} for every $t \in (0,T)$. Moreover, any optimal control $\al$ for \eqref{15} satisfies $\al(u,\omega)=\hat\gamma(u,X(u,\omega))$, $\la_T^t\otimes\PP$ almost surely (a.s.), where $\lambda_T^t$ denotes the Lebesgue measure on $[0,T-t]$.
\end{proposition}
Using Proposition \ref{lem2}, \cite{BBC2017} proves the existence of a solution of MFG under Assumption \ref{assumption1}.
In order to establish  uniqueness of the solution, we need an additional condition. Fix $\eta_0 \in \clp([0,L])$.
\begin{assumption}\label{assumptionU}
For every $(t,\eta,x,u)\in[0,T]\times\calP([0,L])\times[0,L]\times U$,
\begin{align}\label{19q}
b(t,\eta,x,u)&=b(t,\eta_0,x,u),\quad
f(t,\eta,x,u)=f_0(t,\eta,x)+f_1(t,x,u),\\\label{19pp}
y(t,\eta)&=y(t,\eta_0),\quad
r(t,\eta)=r(t,\eta_0).
\end{align}
Moreover, for every $t\in[0,T]$ and $\eta,\eta'\in\calP([0,L])$, $f_0$ and $g$ satisfy the following monotonicity property
 \begin{align}\notag
 \int_0^L [f_0(t,\eta,x)-f_0(t,\eta',x))]d(\eta-\eta')(x)&\ge 0,\\\notag
 \int_0^L (g(\eta,x)-g(\eta',x))d(\eta-\eta')(x)&\ge 0.
 \end{align}
\end{assumption}
Abusing notation, when Assumption \ref{assumptionU} holds, we will write $b(t,x,u) =b(t,\eta_0,x,u)$, $y(t)=y(t,\eta_0)$, and $r(t)=r(t,\eta_0)$.
The following is one of the main results from \cite{BBC2017}.

\begin{proposition}[Theorem 3.1 in \cite{BBC2017}]\label{thm_fixed}
Under Assumption \ref{assumption1}, there exists a solution of the MFG. If in addition Assumption \ref{assumptionU} holds then there is a unique MFG solution.
\end{proposition}
\subsection{Rate control in queues with strategic servers}\label{queapp}

The MFG described above arises from the heavy traffic analysis of a large queuing system that consists of many symmetric strategic servers that are weakly interacting. Consider a collection of $n$ critically loaded single server queues. Given past information, each server controls the arrival and service rate associated with its own queue. In addition the rates depend on time, the individual queue length, and the empirical measures of all the queue states. The servers aim to minimize individual costs, that in particular account for the scaled idleness and rejection processes. The cost also depends on the individual queue state, the control action and the state of the overall system given through the empirical measure of states of all queues.  The main goal is to find asymptotic Nash equilibrium in this game as the system approaches criticality (i.e. heavy traffic limit) and the number of queues approach $\infty$, {\em simultaneously}.  It is shown in \cite{BBC2017} that given a solution of the MFG and an optimal control associated with it of the form in Section \ref{sec2a}, one can construct an asymptotic (in number of players and in heavy traffic limit) Nash equilibrium for the $n$-player game  such that the solution of the MFG and its associated value function approximate the empirical distribution of the states of the queues and the value function of each server. The current work provides a numerical approximation for the solution of the MFG that is needed in order for constructing the above $n$-player asymptotic Nash equilibrium.

\section{Numerical scheme for the MFG}\label{sec3}
In this section we will use the Markov chain approximation method (\cite{Kushner1992}) to construct numerical solutions of the MFG. The main result of the paper
Theorem \ref{thm:thm52} is given here. The numerical scheme is composed of two main steps. First, in Section \ref{sec31}, given a probability measure in $\clp_{T,L}$,
 we construct a finite state, discrete time, controlled Markov chain and provide a numerical scheme to construct a measure over $\calC([0,T]:[0,L])$. Then, in Section \ref{sec32} we show that, under assumptions that include the existence of a unique solution of the MFG, the measure  constructed from the chain converges to the solution of the MFG over a small time interval.
Throughout the section we assume that Assumption \ref{assumption1} is satisfied and that \eqref{19pp} holds. 
Note that we do not assume \eqref{19q} or the monotonicity condition in Assumption \ref{assumptionU}, however
we will introduce additional assumptions as needed.

We now introduce the  controlled Markov chain constructed on some probability space
$(\Om, \clf, \PP)$ that will be used to approximate the solution of the MFG.

\subsection{Approximating controlled Markov chains}\label{sec31}

Fix a discretization parameter $h>0$ such that $L$ is an integer multiple of $h$. Denote the $h$-grid $\{-h,0,h,\dots,L+h\}$ by $\bbS^h$. This is a discretized version of the state space $[0,L]$. Since $0$ and $L$
are reflecting barriers for the state process $X$, we will consider two
 types of transition steps for the approximating chain. The first, which occurs when the chain is away from the boundary,
will be referred to as the {\it rate control step} and the second occurs at the end points $L+h$ and $-h$ and is referred to as the {\it reflection step}.\skp

\noi{\bf Rate control step.} For every $t\in [0,T]$, $\eta\in\calP([0,L])$, $u\in U$, and $x\in\bbS^h_0\doteq\bbS^h\setminus \{-h,L+h\}$ let
\begin{align}\label{equ003}
q^h(t,\eta,u; x,x\pm h)&\doteq \frac{\pm hb(t,\eta,x,u)+\sigma^2}{2\sigma^2},\\\notag
\end{align}
Note that
$$\sum_{y \in \{x\pm h\}} q^h(t,\eta,u; x,y) = 1$$
and that for $0<h<\sigma^2/c_B$, the transition probabilities are positive. Hereafter, these inequalities on $h$ are in force.
 Also, define
\begin{align}\label{eq:sptimsc}
\Delta^h\doteq\frac{h^2}{\sigma^2}.
\end{align}
This will be used to define the continuous time interpolation of the controlled Markov chain. Denote the $\Delta^h$-grid $\{0,\Delta^h,2\Delta^h,\ldots,T-\Delta^h\}$ by $\mathbb{T}^h$.
Note that
$\Delta^h\to 0$ as $h\to 0$.
One can verify that
the following {\it local consistency} conditions (cf.~\cite{Kushner1992}) hold for every $x\in\bbS^h_0$,
\begin{align}\label{num4}
m^h_0(t,\eta,u,x)&\doteq\int_{\bbS^h}(\tilde x-x)q^h(t,\eta,u; x,d\tilde x) = b(t,\eta,x,u)\Delta^h,\\\label{num5}
(\sigma^h_0)^2(t,\eta,u,x)&\doteq\int_{\bbS^h}(\tilde x-x-m^h_0)^2q^h(t,\eta,u; x,d\tilde x) =\sigma^2\Delta^h-(b(t,\eta,u,x)\Delta^h)^2.
\end{align}

\noi{\bf Reflection step.} Such a step occurs only when $x\in\{-h,L+h\}$.
For every $t\in[0,T]$, $\eta\in\calP([0,L])$, and $u\in U$ let
\begin{align}\notag
q^h(t,\eta,u; L+h,L)=q^h(t,\eta,u; -h,0)=1.
\end{align}

We will now define a controlled Markov chain $\{X_n^{h,\nu}\}_{n \in \N_0}$ associated with the parameter $h$, a measure $\nu \in \clp_{T,L}$ and an initial condition $x_0\in [0,L]$. 
We choose to work with a deterministic initial state for simplicity of presentation.  The results continue to hold  when the  initial state is random. In that case, in Construction \ref{const:const1},
one needs an additional initialization step at which one takes a random draw $x_0$ from the initial distribution and then sets $x^h_0=\lfloor x_0/h\rfloor h$, as in the construction below.

We will assume that
$ I(h) \doteq T/\Delta^h $ and $L/h$ are integers. 

\begin{remark}
	The numerical scheme that we develop will be based on controlled Markov chains associated with the probability kernel $q^h$. Such controlled Markov chain based schemes are  closely related to explicit finite difference schemes for parabolic PDE. Although not studied in the current work, one can also consider Markov chain approximation schemes that have behavior similar to that of implicit finite difference schemes. One of the important steps in convergence proofs
 of finite difference schemes is the identification of appropriate stability conditions for space-time discretizations. 
For the Markov chain approximation method, the analogue of such stability conditions are the local consistency condition of the form in \eqref{num4}--\eqref{num5} which form the heart of our convergence proof.  These local consistency requirements in particular imply for our scheme the space-time scaling of the form in \eqref{eq:sptimsc}.  
	
\end{remark}

\begin{construction}
	\label{const:const1}
		$\\$
\begin{itemize}
	\item Define $X_0^{h,\nu} = x^h_0=\lfloor x_0/h\rfloor h$, set $t^{h,\nu}_0=0$ and let $\alpha_{-1}^{h,\nu}$ be a fixed element of $U$.
	\item Having defined for $i=0, 1, \ldots, n$ time instants $t_i^{h,\nu}<T$ and random variables $X_i^{h,\nu}, \alpha^{h,\nu}_{i-1}$ with values in $\bbS^h$ and $U$ respectively, let $\clf_i^{h,\nu} \doteq \sigma \{X^{h,\nu}_j, \alpha^{h,\nu}_{j-1}: j = 0, 1, \ldots i\}$.
	\item Choose the control $\alpha_n^{h,\nu}$ for the $n$-th step that is a $U$ valued $\clf_n^{h,\nu}$ measurable random variable and
	let $X^{h,\nu}_{n+1}$ be such that its conditional distribution 
	given $\clf^{h,\nu}_n$ is
	$q^{h}(t_n^{h,\nu},\nu(t_n^{h,\nu}), \alpha_n^{h,\nu},  X_n^{h,\nu}, \cdot)$, where $\nu(t)$ denotes the marginal distribution of $\nu$
	at time instant $t$. Also define
	\begin{align}\label{eq:tn}
	t^{h,\nu}_{n+1} \doteq t_n^{h,\nu} + \Delta^h 1_{\{X_{n}^{h,\nu} \notin \{-h, L+h\}\}}
	\end{align}
	where the indicator in the above definition will ensure that when we do a continuous time interpolation of the chain, reflection steps `occur instantaneously'. Note that the choice of $\alpha^{h,\nu}_{n}$ is irrelevant if
	$X^{h,\nu}_{n} \in \{-h, L+h\}$.
\end{itemize}
\end{construction}
If $\alpha_n^{h,\nu} = \vartheta(t_n^{h,\nu}, X_{n}^{h,\nu})$ for some $\vartheta: \mathbb{T}^h \times \mathbb{S}^h \to U$ 
then the function $\vartheta$ is referred to as a feedback control.\skp

\noi{\bf Some auxiliary processes.} We will now introduce some processes that will be useful in the analysis of the $h$-th Markov chain. 
Consider the piecewise constant processes
\begin{align} (X^{h,\nu}(t), \alpha^{h,\nu}(t)) \doteq (X^{h,\nu}_{n^{h,\nu}(t)}, \alpha^{h,\nu}_{n^{h,\nu}(t)}),\quad t \in [0,T],\label{eq:eq905z}
\end{align}
where
\begin{equation}
	n^{h,\nu}(t) \doteq 
	 \max\{n: t^{h,\nu}_n = j \Delta^h\}, \; t \in [j\Delta^h, (j+1)\Delta^h),\quad j = 0, \ldots I(h)-1.\label{eq:eq906}\end{equation}
Let $F^{h,\nu}(0)=B^{h,\nu}(0)=Y^{h,\nu}(0)=R^{h,\nu}(0)=0$ and for every $t\in[0,T]$ let
\begin{align}\notag
F^{h,\nu}(t)&\doteq\sum_{j=0}^{n^{h,\nu}(t)-1}\E\left[X^{h,\nu}_{j+1}-X^{h,\nu}_j\mid\calF^{h,\nu}_j\right]1_{\{X^{h,\nu}_j\notin\{-h,L+h\}\}},\\\label{eq:eq907}
B^{h,\nu}(t)&\doteq\frac{1}{\sigma}\sum_{j=0}^{n^{h,\nu}(t)-1}\Big(X^{h,\nu}_{j+1}-X^{h,\nu}_j-\E[X^{h,\nu}_{j+1}-X^{h,\nu}_j\mid\calF^{h,\nu}_j]\Big)1_{\{X^{h,\nu}_j\notin\{-h,L+h\}\}},\\\notag
Y^{h,\nu}(t)&\doteq\sum_{j=0}^{n^{h,\nu}(t)-1}(X^{h,\nu}_{j+1}-X^{h,\nu}_j)1_{\{X^{h,\nu}_j=-h\}} = h\sum_{j=0}^{n^{h,\nu}(t)-1}1_{\{X^{h,\nu}_j=-h\}},\\\notag
R^{h,\nu}(t)&\doteq\sum_{j=0}^{n^{h,\nu}(t)-1}(X^{h,\nu}_{j}-X^{h,\nu}_{j+1})1_{\{X^{h,\nu}_j=L+h\}} = h\sum_{j=0}^{n^{h,\nu}(t)-1}1_{\{X^{h,\nu}_j=L+h\}}.
\end{align}
One can verify that the following representation holds
\begin{align}\label{num11}
(X^{h,\nu},Y^{h,\nu},R^{h,\nu})(t)=\Gamma\left(x^h_0+F^{h,\nu}(\cdot)+\sigma B^{h,\nu}(\cdot)\right)(t),\quad t\in[0,T].
\end{align}
Also, from \eqref{num4} and \eqref{num5} it follows that, on the set $\{X^{h,\nu}_n\notin\{-h,L+h\}\}$ 
\begin{align}\notag
&\E[X^{h,\nu}_{n+1}-X^{h,\nu}_n\mid\calF^{h,\nu}_n]=b(t^{h,\nu}_n,\nu(t^{h,\nu}_n),X^{h,\nu}(t^{h,\nu}_n),\al^{h,\nu}(t^{h,\nu}_n))\Delta^h
,\\\notag
&\E\left[\left[X^{h,\nu}_{n+1}-X^{h,\nu}_n-\E[X^{h,\nu}_{n+1}-X^{h,\nu}_n\mid\calF^{h,\nu}_n] \right]^2\mid \calF^{h,\nu}_n \right]=
\sigma^2\Delta^h+o(\Delta^h).
\end{align}
Define $\clg^{h,\nu}_t \doteq \clf^{h,\nu}_{n^{h,\nu}(t)}$. Then since $n^{h,\nu}(t)$ for each fixed $t$ is a $\{\clf^{h,\nu}_j\}$ stopping time, we have by optional sampling theorem that $B^{h,\nu}(\cdot)$
is a $\{\clg^{h,\nu}_t\}$ martingale. Also, from the above,
\begin{align}\label{num14}
F^{h,\nu}(t)&=\int_0^t b(l^h(s),\nu(l^h(s)),X^{h,\nu}(s),\al^{h,\nu}(s))ds,
\end{align}
where
\begin{align}\notag
l^h(s)\doteq \lfloor s/\Delta^h\rfloor\Delta^h,\qquad s\in[0,T].
\end{align}

\noi{\bf Cost function for the MDP.}
For 
every $(t,x)\in\mathbb{T}^h\times\mathbb{S}^h$ and any admissible control $\alpha^{h,\nu}$ used to construct the $h$-th controlled Markov chain, define the associated cost
\begin{align}\label{num16}
J^{h,\nu}(t,x,\al^{h,\nu})& \doteq\E\left[\int_t^T f(l^h(s),\nu(l^h(s)),X^{h,\nu}(s),\al^{h,\nu}(s))ds+g(\nu(T),X^{h,\nu}(T))\right.\\\notag
&\qquad\left.+\int_t^Ty(s)dY^{h,\nu}(s)+\int_t^Tr(s)dR^{h,\nu}(s)\mid X^{h,\nu}(t)=x\right].
\end{align}
The value function associated with the above cost is given by,
\begin{align}\notag
V^h_{\nu}(t,x)\doteq \inf_{\al}J^{h,\nu}(t,x,\al),
\end{align}
where the infimum is taken over all admissible controls.

We now provide properties of the value function $V^h_\nu$ and the optimal strategy in the $h$-th  MDP. 
For every $(t,x,\nu)\in\mathbb{T}^h\times\mathbb{S}^h_0\times\calPTL$, define the $h$-th finite difference of the value function w.r.t.~$x$, as follows
\begin{align}\label{1009}
\calD^{h}_xV^{h}_\nu(t,x)&\doteq \frac{1}{2h}(V^{h}_\nu(t+\Delta^h,x+h)-V^{h}_\nu(t+\Delta^h,x-h)),
\end{align}
where 
\begin{align}\label{1008}
V^{h}_\nu(t+\Delta^h,L+h)&\doteq r(t+\Delta^h)h+V^{h}_\nu(t+\Delta^h,L),\\\notag
 V^{h}_\nu(t+\Delta^h,-h)&\doteq y(t+\Delta^h)h+V^{h}_\nu(t+\Delta^h,0).
\end{align}

\begin{lemma}\label{lem_52}
The optimal control in the $h$-th  MDP is given in state feedback form as
\begin{align}\label{1006da}
\vartheta^{h,\nu}(t,x)= \hat\al(t,\nu(t),x,\calD^{h}_xV^{h}_\nu(t,x)),\quad (t,x)\in\mathbb{T}^h\times\mathbb{S}^h.
 \end{align}
Letting
 \begin{align}\label{1006d}
\hat\al^{h,\nu}(t,x)\doteq \vartheta^{h,\nu}(l^h(t),x),\quad (t,x)\in [0,T]\times\mathbb{S}^h,
 \end{align}
 there exists a constant $c_d(T)\in(0,\iy)$, such that for every $(t,\nu)\in\mathbb{T}^h\times\calPTL$ and for every $h$, one has,
 \begin{align}\label{1006b}
&V^{h}_\nu(t,X^{h,\nu}(t))+ \sigma\sum_{s\in\mathbb{T}^h, s\ge t}\calD^{h}_xV^{h}_\nu(s,X^{h,\nu}(s))(B^{h,\nu}(s+\Delta^h)-B^{h,\nu}(s))\\\notag
&\quad=g(\nu(T),X^{h,\nu}(T))+\int_t^Tf(l^h(s),\nu(l^h(s)),X^{h,\nu}(s),\hat \al^{h,\nu}(s,X^{h,\nu}(s)))ds\\\notag
&\qquad+\int_t^Ty(s)dY^{h,\nu}(s)+\int_t^Tr(s)dR^{h,\nu}(s),
\end{align}
where $(X^{h,\nu},B^{h,\nu},Y^{h,\nu},R^{h,\nu})$ are as in \eqref{eq:eq905z}-\eqref{eq:eq907} with $\{\alpha^{h,\nu}_n\}$ replaced with the optimal feedback control
$\vartheta^{h,\nu}$ and for all $(t,x)\in\mathbb{T}^h\times\mathbb{S}^h$
\begin{align}\label{1006c}
|\calD^{h}_xV^{h}_\nu(t,X^{h,\nu}(t))|\le c_d(T).
 \end{align}
\end{lemma}
We note that Lemma \ref{lem_52} gives, in an explicit form, the finite difference scheme associated with the dynamic programming equation for the cost function \eqref{num16}. Indeed, recall that the function $\hat \al$ is given in \eqref{19a}, that the gradient $\calD^{h}_xV^{h}_\nu(t,x)$ from \eqref{1009} is calculated based on the values of the value function at time $t+\Delta^h$, and that the integrals in \eqref{1006b} can be written as finite sums over $s\in\mathbb{T}^h\cap[t,T]$. Now, by taking (conditional) expected values in \eqref{1006b}, the finite difference scheme follows from a backwards induction.  
The proof of the lemma is given in Section \ref{sec4}.\skp

\noi{\bf The induced measure $\boldsymbol{\Phi^h(\nu)}$.} Recall from \eqref{eq:eq906} that for $j = 0, 1, \ldots I(h)$, $n^{h,\nu}(j \Delta^h) = \max\{i: t_i^{h,\nu} = j \Delta^h\}$. Let $\{\hat X^{h,\nu}(t)\}_{t\in [0,T]}$ be the continuous stochastic process which is linear on $[j\Delta^h, (j+1)\Delta^h]$
	and equals $X^{h,\nu}_{n^{h,\nu}(j \Delta^h)}$ at $t = j\Delta^h$, for $j=0, \ldots I(h)-1$, where $\{X^{h,\nu}_n\}$ is the controlled Markov chain constructed using the
	optimal feedback control $\hat \al^{h,\nu}$.
	Let
	\begin{align}\label{1000z}
	\Phi^h(\nu)\doteq \PP\circ(\hat X^{h,\nu})^{-1}.
	\end{align}
	Next, we show that,  under suitable conditions,  $\Phi^h$ is a contraction up to an $\clo(h^2)$ term, over a small time interval.  This `almost-contraction' property lies at the heart of our main result, Theorem \ref{thm:thm52}.

\noindent{\bf An almost-contraction property.}
Recall that we assume  that Assumption \ref{assumption1} is satisfied.
In addition, we will make the following assumption on a Lipschitz property of the  function $\hat \al$ from \eqref{19a}.
 \begin{assumption}\label{assu:lips_al}$\,$
There exists $c_\al \in (0,\infty)$ such that for every $t\in[0,T]$, $\eta,\eta'\in\calP([0,L])$, $x,x'\in[0,L]$, and $p,p'\in \R$,
\begin{align}\label{1001}
|\hat\al (t,\eta,x,p)-\hat\al (t,\eta',x',p')|\le c_\al(W_1(\eta,\eta')+|x-x'|+|p-p'|).
\end{align}
\end{assumption}

The following lemma gives a sufficient condition for Assumption \ref{assu:lips_al} to hold.
\begin{lemma}\label{lem_51}
Suppose that the drift and the cost functions satisfy the following properties.
\begin{enumerate}
\item[(a)] For every $(t,\eta,x,\al)\in[0,T]\times\calP([0,L])\times[0,L]\times U$,
\begin{align}\label{1002}
b(t,\eta,x,\al)=b_1(t,\eta,x)+b_2(t)\al.
\end{align}
\item[(b)] There exists $c_{m} \in (0,\infty)$ such that for every $t\in[0,T]$, $\eta,\eta'\in\calP([0,L])$, $x,x'\in[0,L]$, and $\al,\al'\in U$,
\begin{align}\label{1003}
f(t,\eta,x,\al')-f(t,\eta,x,\al)-(\al'-\al)f_\al(t,\eta,x,\al)\ge c_{m}(\al'-\al)^2.
\end{align}
\item[(c)]
The map $\alpha \mapsto f(t,\eta,x,\al,p)$ is continuously differentiable for every $(t,\eta,x,p) \in [0,T]\times\calP([0,L])\times[0,L]\times \R$ and
there exists $c_l \in (0,\infty)$ such that for every $t\in[0,T]$, $\eta,\eta'\in\calP([0,L])$, $x,x'\in[0,L]$, $\al\in U$, and $p,p'\in \R$,
\begin{align}\label{1004}
|f_\al(t,\eta,x,\al,p)-f_\al(t,\eta',x',\al,p')|\le c_l(W(\eta,\eta')+|x-x'|+|p-p'|).
\end{align}
\end{enumerate}
Then Assumption \ref{assu:lips_al} is satisfied.
\end{lemma}
The proof of the lemma is deferred to Section \ref{sec4}.

\begin{remark}\label{rem_new}
The conditions in the above lemma are not new to the literature of MFG. For example, Assumptions (A.1), (A.2), and (A.3) in \cite{Carmona2013} are stronger. Also, parts (b) and (c) above that concern the running cost are imposed by \cite{Gomes2013}, which studies a rate control problem (part (a) is irrelevant for that model).
A basic example that satisfies parts (a)--(c) in the Lemma, in addition to Assumptions \ref{assumption1} and \ref{assumptionU}, is the following
\begin{align}\notag
b(t,\eta,x,\al)&=b_1(t,x)+b_2(t)\al,\\\notag
f(t,\eta,x,\al)&=a_1(t,x)+a_2(t,x)k(\al)+a_3(t)(c_1+a_4(x))\int_0^La_4(y)d\eta(y),\\\notag
g(\eta,x)&=(c_2+a_5(x))\int_0^La_5(y)d\eta(y),\\\notag
y(t,\eta)&=a_6(t),\qquad r(t,\eta)=a_7(t),
\end{align}
where $b_1,a_1,a_2:[0,T]\times[0,L]\to\R$,  $b_2,a_3,a_6, a_7:[0,T]\to\R$, $a_4,a_5:[0,L]\to\R$ are Lipschitz functions,  $c_1,c_2\in\R$,  $k: U \to \R$ is a $C^2$-strictly convex function (e.g.
$k(\alpha) = (\alpha-\alpha_0)^2$ for some $\alpha_0 \in \R$),
and
$a_2(t,x)\ge c_{m}>0$ for all $(t,x)\in [0,T]\times [0,L]$.
We note that, although Assumption \ref{assumptionU} is not explicitly imposed in the current work, we will assume later in the section that the MFG has a unique solution (see Assumption \ref{assu:uniqmfg}) which from Proposition \ref{thm_fixed} holds under Assumptions
\ref{assumption1} and \ref{assumptionU}. For this reason we presented an example that satisfies all three assumptions (i.e., Assumptions \ref{assumption1}, \ref{assumptionU} and \ref{assu:lips_al}). From a modeling perspective, by choosing positive and nondecreasing $a_3$ and $a_4$ and a positive $a_3$, the system planner penalizes all servers collectively for congestion when the empirical measure has high
$a_3$ and $a_4$-moments and in addition it penalizes  individual servers for long queues. Also, when $a_7>0$, rejections of jobs by an individual server are disincentivized and when $a_6<0$, idleness is being rewarded. Finally a convex nondecreasing $k$ assigns costs for increasing the rates.  
\end{remark}


The next result plays an important role in the proof of Theorem \ref{thm:thm52}. Recall Assumption \ref{assumption1}
and \eqref{19pp} are in force.
\begin{proposition}\label{prop_h_cont}
Suppose that Assumption  \ref{assu:lips_al} is satisfied. Then there exist $\hat T>0$, $\hat h>0$ and $q\in(0,1)$, such that for every $T\le\hat T$, $\nu,\nu'\in\calPTL$ and $h\in(0,\hat h\wedge \hat T)$,
\begin{align}\notag
W^2_1(\Phi^h(\nu),\Phi^h(\nu')) \le q \Big(h^2+ W^2_1(\nu,\nu')\Big).
\end{align}
\end{proposition}
The proof of the proposition is given in Section \ref{sec4}.

\subsection{ Approximating the solution of the MFG}\label{sec32}
We now provide the numerical scheme that approximates the solution of the MFG.
 \begin{construction}
	\label{const:const2}
	Let $\hat T, \hat h$ be as in Proposition \ref{prop_h_cont}.
Fix $T<\hat T$ and $(x^h_0,\nu^1,h)\in\mathbb{S}^h\times\calPTL\times (0,\hat h\wedge \hat T)$. Let $\{X^{h,\nu^1}_n\}$ be the $h$-th Markov chain from Construction \ref{const:const1} associated with the optimal control $\hat\al^{h,\nu^1}$. Having defined for $m\in\N$ the process $\{X^{h,\nu^m}_n\}$
, set $\nu^{m+1}\doteq\Phi^h(\nu^m)$ and let $\{X^{h,\nu^{m+1}}_n\}$ be the $h$-th Markov chain from Construction \ref{const:const1} associated with the optimal control $\hat\al^{h,\nu^{m+1}}$.
\end{construction}

With $q\in (0,1)$ as in Proposition \ref{prop_h_cont} we get that for every $h$ as in Construction \ref{const:const2} and every $k\in\N$,
\begin{align}\notag
W^2_1(\Phi^h(\nu^k),\nu^k)=W^2_1(\Phi^h(\nu^k),\Phi^h(\nu^{k-1}))\le q(h^2+W^2_1(\nu^k,\nu^{k-1})).
\end{align}
By iterating this bound we obtain
\begin{align}\label{1060}
W^2_1(\Phi^h(\nu^k),\nu^k)\le\frac{q}{1-q}h^2+q^{k-1}W^2_1(\nu^2,\nu^1)).
\end{align}
Set
\begin{align}\notag
k_h\doteq \min\left\{k\in\N : W^2_1(\Phi^h(\nu^k),\nu^k)\le \frac{2q}{1-q}h^2\right\}
\end{align}
and
\begin{align}\label{eq:eq444}
\nu_h\doteq \nu^{k_h}.
\end{align}
We note that $k_h$ depends also on $\nu^1$, however, it plays no role in the sequel and is therefore omitted from the notation. Processes $(X^{h,\nu_h},Y^{h,\nu_h},R^{h,\nu_h},B^{h,\nu_h})$ are defined as in Construction \ref{const:const1},
replacing $\nu$ with $\nu_h$ and $\alpha_n^{h,\nu}$ with $\hat\alpha_n^{h,\nu_h}$.

As an immediate consequence of the definition of $\nu_h$, we get the following proposition, which is key to the proof of the approximation result in Theorem \ref{thm:thm52} below.
\begin{proposition}\label{prop52}
Suppose that Assumption  \ref{assu:lips_al} is satisfied. Then with $\hat T$ as in Proposition \ref{prop_h_cont},
 for every $T\le\hat T$
\begin{align}\label{1061}
\lim_{h\to 0}\,W^2_1(\Phi^h(\nu_h),\nu_h)=0.
\end{align}
\end{proposition}
For the main result of this section (Theorem \ref{thm:thm52}), in addition to Assumptions \ref{assumption1}, \ref{assu:lips_al} and the property in \eqref{19pp}, we also need the following assumption.
\begin{assumption}
	\label{assu:uniqmfg} 
	There is a unique $\bar \nu \in \clp_{T,L}$ that solves the MFG with initial condition $x$.
\end{assumption}

In order to formalize the main result we introduce the notion of {\it relaxed controls}. The reason is that we need to argue the tightness of control sequences in an appropriate space.
For this, we borrow a relaxed control formulation from \cite[Section 4.3]{BBC2017}. Consider the relaxation of the stochastic control problem in \eqref{13}--\eqref{15} where the control space $U$ is replaced by
$\clp(U)$, the drift function $b$ is replaced by the function $b_{\clr}:[0,T]\times [0,L]\times \clp(U) \to \R$
defined as
$$b_{\clr}(t,x, r) \doteq \int_U b(t, x, u) r(du),
,$$
and the running cost $f$ is replaced by $f_{\clr}: [0,T]\times \clp([0,L])\times [0,L]\times \clp(U) \to \R$, defined as
$$f_{\clr}(t,\eta, x, r) \doteq \int_U f(t, \eta, x, u) r(du)
.$$
Finally, we replace the class of admissible controls $\cla(\Xi, t,x, \bar \nu)$ by $\cla_{\clr}(\Xi, t,x, \bar \nu)$ of pairs $(\alpha_{\clr}, Z)$ that are similar to pairs 
$(\alpha, Z)$ introduced above \eqref{13} except that $\alpha_{\clr}$ is $\clp(U)$ valued rather than $U$ valued and in \eqref{13} we replace
$\bar b(u) = b(u, X(u), \alpha(u))$ with $b_{\clr}(u, X(u), \alpha_{\clr}(u))$. The corresponding cost function
$J_{\bar \nu, \clr}$ is defined by \eqref{eq:eq234} with $f$ replaced by $f_{\clr}$. The value function in this relaxed formulation, denoted as $V_{\bar \nu, \clr}$,
is given by \eqref{15} with $\cla$ replaced by $\cla_{\clr}$. Define the function $h_{\clr}$ by \eqref{19}, replacing $(f,b)$ with $(f_{\clr},b_{\clr})$.
Then, from Assumption \ref{assumption1}(b),
$$H(t,\eta,x,p) = \inf_{u \in U} h(t,\eta, x, u,p) = \inf_{r \in \clp(U)} h_{\clr}(t,\eta, x, r,p).$$
To see the last equality note that
$$\inf_{r \in \clp(U)} h_{\clr}(t,\eta, x, r,p)=\inf_{r\in\calP(U)}\int h(t,\eta, x, u,p)dr(u)\ge \inf_{u\in U}h(t,\eta, x, u,p)$$ and on the other hand 
$$ \inf_{u\in U}h(t,\eta, x, u,p)= \inf_{u\in U}\int h(t,\eta, x, a,p)d\delta_u(a)\ge \inf_{r \in \clp(U)} h_{\clr}(t,\eta, x, r,p).$$

Therefore, $V_{\nu}$ and $V_{\nu, \clr}$ are both solutions of the partial differential equation \eqref{HJB1}-\eqref{HJB2}. In view of the uniqueness result given in Proposition \ref{lem2}, $V_{\nu} = V_{\nu, \clr}$.

\begin{remark}
	\label{rem:rem629}
Recall from Proposition \ref{thm_fixed} that 
Assumption \ref{assu:uniqmfg} is satisfied if in addition to Assumption \ref{assumption1},
Assumption \ref{assumptionU} holds.
Also, from Assumptions \ref{assumption1} and \ref{assu:uniqmfg} it follows from arguments as in the proof of Proposition 3.1 in \cite{BBC2017} (see also the statement of Proposition \ref{lem2} here),
that there is a continuous map 
$\gamma:[0,T] \times [0,L] \to U$ such that if there exist a system $\Xi$ and an $(\alpha_{\clr},Z) \in \cla_{\clr}(\Xi,0,x,\nu)$ such that
\begin{equation}\label{eq:eq637c}
 Z(t)\equiv (X,Y,R)(t)=\Gamma\left(x+\int_0^\cdot
b_{\clr}(t, \nu(t), X(t), \alpha_{\clr}(t))dt +\sigma B(\cdot) \right)(t),\; t \in [0,T],
\end{equation}
 $\PP \circ X^{-1} =  \nu$, and
\begin{equation}\label{eq:eq716b}
	V_{\nu}(0,x) = J_{\nu}(0, x, \alpha_{\clr},Z).
\end{equation}
then $\alpha_{\clr}(t, \omega) = \delta_{\gamma(t, X(t, \omega))}$, $\lambda_T^0\otimes \PP$-a.s., where recall that $\lambda_T^0$ denotes the Lebesgue measure on $[0,T]$, and $\nu = \bar \nu$.
\end{remark}

%
%
%
%
%
We now present the main result of the paper. 
Let $\clm(U \times [0,T])$ be the space of finite measures on $U \times [0,T]$ equipped with the topology of weak convergence. 
Define $\clm(U\times [0,T])$ valued random variable $\hat m^{h,\nu_h}$ as
$$\hat m^{h,\nu_h}(du\, ds) \doteq \delta_{\hat \alpha^{h,\nu_h}(s,X^{h,\nu_h}(s))}(du) ds.$$  Also, recall Assumption \ref{assumption1}
and \eqref{19pp} are in force throughout this section.
\begin{theorem}
	\label{thm:thm52}
Suppose that $T\le \hat T$ where $\hat T$ is as in Proposition \ref{prop_h_cont}.
Also suppose that  Assumptions
\ref{assu:lips_al} and \ref{assu:uniqmfg} are satisfied. Recall the processes $(X^{h,\nu_h}, Y^{h,\nu_h}, R^{h,\nu_h}, B^{h,\nu_h})$ introduced below Construction \ref{const:const2} and
consider a sequence $\{h\}\to 0$.
Then
the sequence $$(X^{h,\nu_h}, Y^{h,\nu_h}, R^{h,\nu_h}, B^{h,\nu_h}, \hat m^{h,\nu_h}),$$  converges in distribution to
$(X,Y, R, B, m)$ in $\cld([0,T]:\R^4)\times \clm(U \times [0,T])$ and
\begin{align}\label{1070}
\lim_{h\to0}\nu_h= \nu
\end{align}
 in $\clp_{T,L}$ where the
limit processes defined on some probability space $(\Omega, \clf, \PP)$ satisfy the following.
\begin{description}
	\item[(a)] $B$ is a $\clg_t \doteq \sigma\{ B(s), X(s), Y(s), R(s), m(A\times[0,s]): s \le t, A \in \clb(U)\}$ Brownian motion and so
	 $\Xi \doteq (\Omega, \clf, \{\clg_t\}, \PP, B)$ is a system.
	\item[(b)] Disintegrating $m(du\, ds) = m_s(du) ds$, the following relationship holds a.s.
	$$Z(t)\doteq (X(t), Y(t), R(t)) = \Gamma\left(x + \int_0^{\cdot}b_{\clr}(s, \nu(s), X(s), m_s) ds + \sigma B(\cdot)\right)(t), \; t \in [0,T].$$
	\item[(c)] $\PP\circ (X)^{-1} = \nu$.
	\item[(d)] The pair $(m,Z) \in \cla_{\clr}(\Xi,0,x,\nu)$ and $V_{\nu}(0,x) = J_{\nu, \clr}(0, x, m,Z)$. In particular, with $\gamma$ as in Remark \ref{rem:rem629},
	$m(du\, ds) = \delta_{\gamma(s, X(s))}(du) ds$, and $\nu = \bar \nu$, the unique solution of the MFG.
\end{description}
%
%
%
 \end{theorem}
\noi{\bf Proof.}
Since many steps in the proof are quite standard we will only provide details where appropriate. Using properties
\eqref{num4} and \eqref{num5} of the controlled transition probability kernel it can be argued (cf.~Proof of \cite[Theorem 9.4.1]{Kushner1992}) that $\{(F^{h,\nu_h}, B^{h,\nu_h})\}_{h>0}$ are tight in $\cld([0,T]:\R^2)$. Using this tightness property along with the continuity of the Skorohod map (Lemma \ref{lem_Skorohod}) it now follows that $\{(X^{h,\nu_h}, Y^{h,\nu_h}, R^{h,\nu_h}, B^{h,\nu_h})\}_{h>0}$ is tight in
$\cld([0,T]:\R^4)$. In fact, this sequence is $\clc$-tight.
Next note that, since $\hat m^{h,\nu_h}(U\times [0,T]) =T$ and $U$ is compact, the sequence $\{\hat m^{h,\nu_h}\}_{h>0}$
is tight in $\clm(U\times [0,T])$. Also, recalling the definition of the interpolated processes (right before \eqref{1000z}) it can be checked that
\begin{equation}\label{eq:eq539}
|X^{h,\nu_h} - \hat X^{h,\nu_h}|_T \to 0, \mbox{ in probability as } h\to 0.\end{equation}
Combining this with the tightness of $\{X^{h,\nu_h}\}$ and the fact that $\Phi^h(\nu_h) = \PP \circ (\hat X^{h,\nu_h})^{-1}$ gives the relative compactness of $\{\Phi^h(\nu_h)\}$ in $\clp_{T,L}$. Suppose now that along a subsequence (relabeled again as $\{h\}$)
\begin{equation}
	\label{eq:eq1002}(X^{h,\nu_h}, Y^{h,\nu_h}, R^{h,\nu_h}, B^{h,\nu_h}, \hat m^{h,\nu_h}) \Rightarrow (X,Y, R, B, m)\end{equation}
	and
	\begin{equation}\label{eq:eq1002b}
	 \Phi^h(\nu_h) \to \nu.
	\end{equation}
Then by \eqref{1061} we also have that
\begin{align}\label{1071}
\nu_h\to \nu.
\end{align}
By standard martingale methods it follows that $B$ is a $\{\clg_t\}$ Brownian motion (see e.g.~Proof of \cite[Theorem 9.4.1]{Kushner1992}) proving part (a) of the theorem.

Using \eqref{num14},
$F^{h,\nu_h}$ converges, along with the above processes, in distribution to
$$\int_{U\times [0,\cdot]} b(s, \nu(s), X(s), u) m(du\, ds) = \int_{0}^{\cdot} b_{\clr}(s, \nu(s), X(s), m_s) ds,$$
where $m(du ds) = m_s(du) ds$. Using the continuity property of the Skorohod map we now get (b).

Also, part (c) is immediate on using \eqref{eq:eq539}, recalling that $\Phi^h(\nu_h)$ is the probability law of $\hat X^{h,\nu_h}$, and by \eqref{1071}.

Clearly $(m,Z) \in \cla_{\clr}(\Xi, 0,x, \nu)$. We will now argue the first statement in (d), namely
\begin{equation}\label{eq:eq543}
	V_{\nu}(0,x) = J_{\nu, \clr}(0, x, m,Z).
\end{equation}
This together with Remark  \ref{rem:rem629} will prove the second statement in (d) for the subsequence.
Since the convergent subsequence was arbitrary, we will get the convergence asserted in the statement of the theorem and complete the proof.

\noi {\bf Proof of \eqref{eq:eq543}.} For the rest of the proof we consider the subsequence along which the convergence in \eqref{eq:eq1002} and \eqref{eq:eq1002b} holds.
Using arguments similar to those in Proposition 4.2 in \cite{BBC2017} it can be checked that
\begin{equation}\label{eq:eq800}
	\lim_{h\to 0} J^{h,\nu_h}(\hat \alpha^{h,\nu_h}) = J_{\nu, \clr}(0, x,m,Z),
\end{equation}
where we have suppressed $(0,x_0^h)$ from the notation in $J^{h,\nu_h}$.
Next, let $(\beta, Z) \in \cla(\bar \Xi, 0, x, \nu)$ for some system $\bar \Xi = (\bar \Omega, \bar \clf, \{\bar \clf_t\}, \bar \PP, \bar B)$.
	We will now show that for every $\eps_0 > 0$ there is a sequence of controls $\{\beta^{h,\nu_h}_i\}$ for the controlled Markov chains defined as in Construction \ref{const:const1} with the sequence $\{\nu_h\}$
, such that
	\begin{equation}
		\label{eq:eq748}
		\limsup_{h\to 0} J^{h,\nu_h}(\beta^{h,\nu_h}) \le J_{\nu}(0, x, \beta, Z) + \eps_0.
	\end{equation}
	Note that from the optimality property of $\hat \alpha^{h,\nu_h}$
	$$\limsup_{h\to 0} J^{h,\nu_h}(\hat\alpha^{h,\nu_h}) \le \limsup_{h\to 0} J^{h,\nu_h}(\beta^{h,\nu_h}).$$
	Combining the above inequality with \eqref{eq:eq800} and \eqref{eq:eq748} we now get that
	$$V_{\nu}(0,x) \le J_{\nu, \clr}(0, x, m,Z) = \lim_{h\to 0} J^{h,\nu_h}(\hat \alpha^{h,\nu_h}) \le \limsup_{h\to 0} J^{h,\nu_h}(\beta^{h,\nu_h}) \le J_{\nu}(0, x, \beta, Z) + \eps_0 .$$
	Since $\eps_0>0$, the system $\bar\Xi$, and $(\beta, Z) \in \cla(\bar \Xi, 0, x, \nu)$ are arbitrary we have \eqref{eq:eq543}, completing the proof of the theorem.
	
	We now prove \eqref{eq:eq748}. Using arguments as in the proof of \cite[Theorem 10.3.1]{Kushner1992}
	it can be shown that there is a $\theta_1:[0,\infty) \to [0,\infty)$ such that $\theta_1(\kappa)\to 0$ as $\kappa \to 0$
	and for every $\eps >0$
	there is a system $\Xi^{\eps}\doteq (\Omega^{\eps},\calF^{\eps},\{\calF_t^{\eps}\},\PP^{\eps}, B^{\eps})$
	and $(\beta^{\eps}, Z^{\eps}) \in \cla(\Xi^{\eps}, 0,x,\nu)$
	with the following properties
	\begin{itemize}
		\item $Z^{\eps}$ satisfies the following equation for $t \in [0,T]$.
		\begin{equation} \label{eq:eq1019} Z^{\eps}(t) \equiv (X^{\eps},Y^{\eps},R^{\eps})(t)=\Gamma\left(x+\int_0^\cdot
		b(t, \nu(s), X^{\eps}(s), \beta^{\eps}(s))ds +\sigma B^{\eps}(\cdot) \right)(t).\end{equation}
		\item For some $\delta > 0$, $\beta^{\eps}$ is piecewise constant on intervals of the form $[l\delta,(l+1)\delta)$, $l=0,1,\ldots, T/\delta$. For some finite set $U^{\eps} \subset U$, $\beta^{\eps}(s)$ takes values in $U^{\eps}$ for every
		$s \in [0,T]$.
		\item For some $\theta>0$, for each $u \in U^{\eps}$
		\begin{align}
		\PP^{\eps}(\beta^\eps(l\delta)=u|B^\eps(s),s\le l\delta, \beta^\eps(j\delta),j<l)&=\PP^{\eps}(\beta^\eps(l\delta)=u|B^\eps(p\theta),p\theta\le l\delta, \beta^\eps(j\delta),j<l)\nonumber\\
		&= F^{\eps}_u(B^\eps(p\theta),p\theta\le l\delta, \beta^\eps(j\delta),j<l),\label{num28}
		\end{align}
		where for suitable $\hat I,\hat L \in \N$, $F^{\eps}_u: \R^{\hat I} \times (U^{\eps})^{\hat L} \to [0,1]$ is a measurable function such that
		$F_u(\cdot, \mathbf{u})$ is continuous on $\R^{\hat I}$ for every $\mathbf{u} \in (U^{\eps})^{\hat L}$.
		\item Letting $m^{\eps}(du dt) \doteq \delta_{\beta^{\eps}(t)}(du) dt$, $m(du dt) \doteq
		\delta_{\beta(t)}(du) dt$, as $\eps \to 0$,
		$$(X^{\eps},Y^{\eps}, R^{\eps}, B^{\eps}, m^{\eps}) \Rightarrow (X,Y, R, B, m)$$
		in $\cld([0,T]:\R^4)\times \clm(U\times [0,T])$.
		\item $J_\nu(0,x,\beta^{\eps},Z^{\eps}) \le J_\nu(0,x,\beta,Z) + \theta_1(\eps)$.
	\end{itemize}

We will now use the piecewise constant control $\beta^\eps$ to construct a collection of control sequences $\{\beta^{h,\nu_h}_i\}$ as stated above \eqref{eq:eq748} and for which
\begin{equation}
	\label{eq:eq914}
	\lim_{h\to 0} J^{h,\nu_h}(\beta^{h,\nu_h}, \nu_h) = J_{\nu}(0, x, \beta^{\eps}, Z^{\eps}).
\end{equation}
Note that since $\theta_1(\eps)\to 0$ as $\eps \to 0$, this will prove \eqref{eq:eq748} and complete the proof of the theorem. The construction is carried out as follows.
\begin{itemize}
	\item Define $X^{h,\nu_h}_0 = x^{h,\nu_h}$, $t_0^1=0$ and let $\beta_{-1}^{h,\nu_h}$ be a fixed element of $U^{\eps}$.
	\item Having defined for $i=0, 1, \ldots, n$ time instants $t_i^{h,\nu_h}<T$ and random variables $X_i^{h,\nu_h}, \beta^{h,\nu_h}_{i-1}$ with values in $\bbS^h$ and $U^{\eps}$ respectively, let $\clf_i^{h,\nu_h} \doteq \sigma \{X^{h,\nu_h}_j, \beta^{h,\nu_h}_{j-1}: j = 0, 1, \ldots i\}$.
	\item Choose the control $\beta_n^{h,\nu_h}$ for the $n$-th step that is a $U^{\eps}$ valued $\clf_n^{h,\nu_h}$ measurable random variable as follows:
	\begin{itemize}
		\item If $t_{n-1}^{h,\nu_h}$ and $t_n^{h,\nu_h}$ both lie in $[j\delta, (j+1)\delta)$ for some $j$, set $\beta^{h,\nu_h}_n = \beta^{h,\nu_h}_{n-1}$.
		\item If $t_{n-1}^{h,\nu_h}< l\delta \le t_n^{h,\nu_h} < (l+1)\delta$ for some $l$, choose $\beta^{h,\nu_h}_n$ according the the conditional distribution
		$$P\left(\beta^{h,\nu_h}_n = u\mid X_i^{h,\nu_h}, \beta^{h,\nu_h}_{i-1}, 0 \le i \le n\right) = F_u^{\eps}(B^{h,\nu_h}(p\theta),p\theta\le l\delta, \beta^{h,\nu_h}(j\delta),j<l),$$
		where $B^{h,\nu_h}$ is defined as in \eqref{eq:eq907}, $\beta^{h,\nu_h}(s) = \beta^{h,\nu_h}_{n^{h,\nu_h}(s)}$, and $n^{h,\nu_h}(\cdot)$ is as in \eqref{eq:eq906}.
	\end{itemize}
	\item Let $X^{h,\nu_h}_{n+1}$ be such that the conditional distribution of $X^{h,\nu_h}_{n+1}$ 
	given $\clf^{h,\nu_h}_n$ equals $\\$
	$q^{h}(t_n^{h,\nu_h}, \nu_h(t_n^{h,\nu_h}), \beta_n^{h,\nu_h}, X_n^{h,\nu_h}, \cdot)$, where $\nu_h$ is as introduced above Theorem \ref{thm:thm52}.
	 Also define
	$$t^{h,\nu_h}_{n+1} \doteq t_n^{h,\nu_h} + \Delta^h 1_{\{X_{n}^{h,\nu_h} \notin \{-h, L+h\}\}}.$$
\end{itemize}
Now define processes $X^{h,\nu_h}, F^{h,\nu_h}, Y^{h,\nu_h}, R^{h,\nu_h}$ as in \eqref{eq:eq905z}--\eqref{eq:eq907}. Exactly as in the first part of the proof we now have that $\{(X^{h,\nu_h}, Y^{h,\nu_h}, R^{h,\nu_h}, B^{h,\nu_h})\}_{h}$ is $\clc$-tight in
$\cld([0,T]:\R^4)$ and the sequence $\{m^{h,\nu_h}\}_{k\ge 1}$, where $m^{h,\nu_h}(du ds) = \delta_{\beta^{h,\nu_h}(s)}(du) ds$
is tight in $\clm(U\times [0,T])$. Arguing as before, if along a further subsequence
the convergence \eqref{eq:eq1002} holds (with $\hat m^{h,\nu_h}$ replaced with $m^{h,\nu_h}$) then parts (a) and (b) as in the statement of Theorem \ref{thm:thm52} are satisfied with $\nu$ as in \eqref{eq:eq1002b}.
Using the continuity property of $F_u$ and the fact that the control is piecewise constant with values in a finite set it follows that
(cf.~Proof of \cite[Theorem 10.5.2]{Kushner1992}) $(B, m)$ has the same distribution as $(B^{\eps}, m^{\eps})$. By unique solvability of \eqref{eq:eq1019}, that follows from the
Lipschitz property of $b$ (Assumption \ref{assumption1}) and the Lipschitz property of the Skorohod map (Lemma \ref{lem_Skorohod}) we now have that
$(Z^{\eps}, m^{\eps})$ has the same law as $(Z, m)$ for every limit point of the chosen further subsequence. Since the chosen further subsequence was arbitrary, this proves the weak convergence of $(Z^{h,\nu_h}, m^{h,\nu_h})$ to $(Z^{\eps},m^{\eps})$ along the subsequence fixed above \eqref{eq:eq800} and
 arguing once again as in the proof of Proposition
4.2 in \cite{BBC2017} we have the convergence of costs
as in \eqref{eq:eq914}, completing the proof of the theorem.
\hfill$\Box$

\section{Proofs of results from Section \ref{sec3}} \label{sec4}
In this section we provide the  proofs of Lemmas \ref{lem_52}, \ref{lem_51} and Proposition \ref{prop_h_cont}.\skp

\noi {\bf Proof of Lemma \ref{lem_52}.} 
We start by analyzing the evolution of the value function $V^{h}_\nu$. Clearly, $V^{h}_\nu(T,x) = g(\nu(T),x)$
for all $x \in \mathbb{S}^h_0$. Using backwards induction, we get that for any 
$(t,x,\nu)\in\mathbb{T}^h\times\mathbb{S}^h_0\times\calPTL$, one has
from the definition of $\hat\al$ (cf.~\eqref{19a}), $q^h$ (see \eqref{equ003}) and \eqref{1006da}-\eqref{1006d},
\begin{align}\label{1007}
V^{h}_\nu(t,x)&=\min_{u\in U}\Big\{\Delta^hf(t,\nu(t),x,u)\\\notag
&\quad\qquad\quad+q^h(t,\nu(t),u;x,x+h)V^{h}_\nu(t+\Delta^h,x+h)\\\notag
&\quad\qquad\quad+q^h(t,\nu(t),u;x,x-h)V^{h}_\nu(t+\Delta^h,x-h)\Big\}\\\notag
&=\Delta^hf(t,\nu(t),x,\hat\al^{h,\nu}(t,x))\\\notag
&\quad+q^h(t,\nu(t),\hat\al^{h,\nu}(t,x);x,x+h)V^{h}_\nu(t+\Delta^h,x+h)\\\notag
&\quad+q^h(t,\nu(t),\hat\al^{h,\nu}(t,x);x,x-h)V^{h}_\nu(t+\Delta^h,x-h),
\end{align}
where $\hat\al^{h,\nu}$ is as in \eqref{1006d}.
The above identity in particular shows that $\vartheta^{h,\nu}$ gives an optimal feedback control.

We now use \eqref{1007} to prove \eqref{1006b}.
Define the following $h$-th finite differences for $(t,x)\in\mathbb{T}^h\times \mathbb{S}^h_0$,
\begin{align}\notag
\calD^{h}_tV^{h}_\nu(t,x)&\doteq \frac{1}{\Delta^h}(V^{h}_\nu(t+\Delta^h,x)-V^{h}_\nu(t,x)),\\\notag
\calD^{h}_{xx}V^{h}_\nu(t,x)&\doteq \frac{1}{h^2}(V^{h}_\nu(t+\Delta^h,x+h)-2V^{h}_\nu(t+\Delta^h,x)+V^{h}_\nu(t+\Delta^h,x-h)).
\end{align}
Simplifying \eqref{1007} by using \eqref{1008}, we get that
 \begin{align}\notag
   \calD^h_tV^{h}_\nu(t,x) =
        -H(t,\nu(t),x,\calD^{h}_x(t,x)) -\frac{1}{2}\sigma^2\calD^h_{xx}V^{h}_\nu(t,x).
  \end{align}
Notice also that $X^{h,\nu}(t+\Delta^h)\ne X^{h,\nu}(t)$ if and only if $R^{h,\nu}(t+\Delta^h)-R^{h,\nu}(t)=Y^{h,\nu}(t+\Delta^h)-Y^{h,\nu}(t)=0$, in which case, one has
\begin{align}\notag
&V^{h}_\nu(t+\Delta^h,X^{h,\nu}(t+\Delta^h))-V^{h}_\nu(t,X^{h,\nu}(t))\\\notag
&\quad=\calD^h_tV^{h}_\nu(t,X^{h,\nu}(t))\Delta^h+\calD^{h}_xV^{h}_\nu(t,X^{h,\nu}(t))(X^{h,\nu}(t+\Delta^h)-X^{h,\nu}(t))+\frac{1}{2}h^2\calD^h_{xx}V^{h}_\nu(t,X^{h,\nu}(t)).
\end{align}
This can be easily verified by considering separately the cases $X^{h,\nu}(t+\Delta^h)=X^{h,\nu}(t)\pm h$.
In case that $R^{h,\nu}(t+\Delta^h)-R^{h,\nu}(t)=h$, that is $X^{h,\nu}(t+\Delta^h)= X^{h,\nu}(t)=L$, \eqref{1008} implies that
\begin{align}\notag
&V^{h}_\nu(t+\Delta^h,X^{h,\nu}(t+\Delta^h))-V^{h}_\nu(t,X^{h,\nu}(t))\\\notag
&\quad =V^{h,\nu}(t+\Delta^h,L)-V^{h,\nu}(t+\Delta^h,L+h)+V^{h,\nu}(t+\Delta^h,L+h)-V^{h,\nu}(t,L)\\\notag
&\quad=-r(t+\Delta^h)h+\calD^h_tV^{h}_\nu(t,X^{h,\nu}(t))\Delta^h+\calD^{h}_xV^{h}_\nu(t,X^{h,\nu}(t))h+\frac{1}{2}h^2\calD^h_{xx}V^{h}_\nu(t,X^{h,\nu}(t))
\end{align}
and finally, in case that $Y^{h,\nu}(t+\Delta^h)-Y^{h,\nu}(t)=h$,
one similarly has
\begin{align}\notag
&V^{h}_\nu(t+\Delta^h,X^{h,\nu}(t+\Delta^h))-V^{h}_\nu(t,X^{h,\nu}(t))\\\notag
\notag
&\quad=-y(t+\Delta^h)h+\calD^h_tV^{h}_\nu(t,X^{h,\nu}(t))\Delta^h+\calD^{h}_xV^{h}_\nu(t,X^{h,\nu}(t))(-h)+\frac{1}{2}h^2\calD^h_{xx}V^{h}_\nu(t,X^{h,\nu}(t)).
\end{align}
From the last two equalities, we get that,
\begin{align}\notag
&V^{h}_\nu(t+\Delta^h,X^{h,\nu}(t+\Delta^h))-V^{h}_\nu(t,X^{h,\nu}(t))\\\notag
&\quad= -H(t,\nu(t),X^{h,\nu}(t),\calD^{h}_xV^{h}_\nu(t,X^{h,\nu}(t)))\Delta^h + \calD^{h}_xV^{h}_\nu(t,X^{h,\nu}(t))(X^{h,\nu}(t+\Delta^h)-X^{h,\nu}(t))\\\notag
&\qquad-(y(t+\Delta^h)+
\calD^{h}_xV^{h}_\nu(t,X^{h,\nu}(t))
)(Y^{h,\nu}(t+\Delta^h)-Y^{h,\nu}(t))\\\notag
&\qquad-(r(t+\Delta^h)-\calD^{h}_xV^{h}_\nu(t,X^{h,\nu}(t)))(R^{h,\nu}(t+\Delta^h)-R^{h,\nu}(t))\\\notag
&\quad= -f(t,\nu(t),X^{h,\nu}(t),\hat \al^{h,\nu}(t,X^{h,\nu}(t)))\Delta^h + \sigma\calD^{h}_xV^{h}_\nu(t,X^{h,\nu}(t))(B^{h,\nu}(t+\Delta^h)-B^{h,\nu}(t))\\\notag
&\qquad-y(t+\Delta^h)(Y^{h,\nu}(t+\Delta^h)-Y^{h,\nu}(t))-r(t+\Delta^h)(R^{h,\nu}(t+\Delta^h)-R^{h,\nu}(t)).
\end{align}
Summing up the terms over $t\in\mathbb{T}^h$, one gets \eqref{1006b}.
We will postpone the proof of \eqref{1006c} to the end of the paper. 
\hfill$\Box$\skp

\noi {\bf Proof of Lemma \ref{lem_51}.}
Fix $t\in[0,T]$, $\eta,\eta'\in\calP([0,L])$, $x,x'\in[0,L]$, and $p,p'\in \R$. Denote $\al=\hat\al(t,\eta,x,p)$ and $\al'=\hat\al(t,\eta',x',p')$. 
Recall the definition of $h$ from \eqref{19}. By \eqref{1002}, \eqref{1003}, and the definition of $\al'$, we get that
\begin{align}\notag
h(t,\eta',x',\al,p')\ge h(t,\eta',x',\al',p')\ge h(t,\eta',x',\al,p')+(\al'-\al)h_\al(t,\eta',x',\al,p')+c_{m}|\al'-\al|^2.
\end{align}
From the minimizing property of $\alpha$ we see  that $(\al'-\al)h_\al(t,\eta,x,\al,p)\ge 0$. Subtracting this term from the right side of the above
\begin{align}\notag
c_{m}|\al'-\al|^2&\le |\al'-\al|\cdot|h_\al(t,\eta',x',\al,p')-h_\al(t,\eta,x,\al,p)|\\\notag
&\le c |\al'-\al|(W_2(\eta,\eta')+|x-x'|+|p-p'|) ,
\end{align}
where $c=c_l+\sup_{t\in[0,T]}|b_2(t)|$ and the second inequality follows by \eqref{1002} and \eqref{1004} . Result follows on dividing both sides by $|\al'-\al|/c_{m}$.
\hfill$\Box$\skp

\noi {\bf Proof of Proposition \ref{prop_h_cont}.}  We begin by introducing a coupling between two optimally controlled chains, one associated with $\nu$ and the other with $\nu'$.

\skp
\noi{\bf Coupling.} Fix $x\in[0,L]$ and $\nu,\nu'\in\calPTL$. Let $\{X^{\nu}_n\}$ and $\{X^{\nu'}_n\}$ be the Markov chains from Construction \ref{const:const1} associated with the parameter $h$ and the optimal strategies given by \eqref{1006d}. Denote by $\Sigma^{\nu}=(n^{\nu},F^{\nu},B^{\nu},X^{\nu},Y^{\nu},R^{\nu})$ and $\Sigma^{\nu'}=(n^{\nu'},F^{\nu'},B^{\nu'},X^{\nu'},Y^{\nu'},R^{\nu'})$ the processes that were defined immediately after Construction \ref{const:const1}, where we suppressed the index $h$ since it is fixed in the rest of the proof.  Also, denote
\begin{align}\notag
b^{\nu}(t)\doteq b(l^h(t),\nu(l^h(t)),X^{\nu}(t),\hat\al^{\nu}(t,X^{\nu}(t))),
\end{align}
where recall that $l^h(t)=\lfloor t/\Delta^h\rfloor \Delta^h$. Similarly define $b^{\nu'}$. We now define a coupling of the chains through a time change of an underlying Markov chain $\{(Z^{\nu}_n,Z^{\nu'}_n)\}$. The main idea in the construction of the latter Markov chain is to keep track of the proper time. Whenever an `instantaneous jump' occurs for only one of the $Z$-processes, the other process has a degenerate step, that is, it remains at the same position. Therefore, we use two sequences of times. The first, which we refer as {\it time instants}, $\{t_n\}$ has the same role as in \eqref{eq:tn}. The second is referred as {\it time steps} and denoted by $\{(N^{\nu}_n,N^{\nu'}_n)\}$. Each of the components counts how many non-degenerate steps the respective $Z$ process has taken so far.

Set $Z^\nu_0=Z^{\nu'}_0=x_0$, $t_0=0$, and $ N^{\nu}_0= N^{\nu'}_0=0$. Having defined for $i=0,1,\ldots,n$ time instants $t^{\nu}_i<T$, time steps $ N^{\nu}_i, N^{\nu'}_i\in\N$, and random variables $Z^\nu_i,Z^{\nu'}_i$ with values in $\mathbb{S}^h$, define them for the $(n+1)$-th step as follows.
\begin{itemize}
	\item If $Z^{\nu}_{n},Z^{\nu'}_{n}\notin \{-h,L+h\}$, then
\begin{align}\notag
 &(Z^{\nu}_{n+1},Z^{\nu'}_{n+1})=(Z^{\nu}_{n},Z^{\nu'}_{n}) +
  \begin{cases}
    (h,h), & \text{w.p. } (h\min\{b^{\nu}(t_n),b^{\nu'}(t_n)\}+\sigma^2)/(2\sigma^2), \\
     (h,-h), & \text{w.p. } h( b^{\nu}(t_n)- b^{\nu'}(t_n))^+/(2\sigma^2), \\
        (-h,h), & \text{w.p. } h(b^{\nu}(t_n)- b^{\nu'}(t_n))^-/(2\sigma^2),\\
        (-h,-h), & \text{w.p. } (-h\max\{b^{\nu}(t_n),b^{\nu'}(t_n)\}+\sigma^2)/(2\sigma^2),
\end{cases}
\end{align}
where w.p. stands for `with probability', and
\begin{align}\notag
t_{n+1}=t_n+\Delta^h,\qquad
 N^{\nu}_{n+1}= N^{\nu}_{n}+1,\qquad&\text{and}\qquad  N^{\nu'}_{n+1}= N^{\nu'}_{n}+1,
\end{align}
where $x^+=\max\{0,x\}$ and $x^-=\max\{0,-x\}$.

\item
If $Z^{\nu}_{n}\notin\{-h,L+h\}$ and $Z^{\nu'}_{n}\in\{-h,L+h\}$ then
\begin{align}\notag
 &(Z^{\nu}_{n+1},Z^{\nu'}_{n+1})=(Z^{\nu}_{n},Z^{\nu'}_{n})+
 \left(0, h(-1)^{\one_{\{Z^{\nu'}_{n}=L+h\}}}\right)
\end{align}
and
\begin{align}\notag
t_{n+1}=t_n,\qquad 
 N^{\nu}_{n+1}= N^{\nu}_{n},\qquad&\text{and}\qquad  N^{\nu'}_{n+1}= N^{\nu'}_{n}+1.
\end{align}
The transition probabilities when $Z^{\nu'}_{n}\notin\{-h,L+h\}$ and $Z^{\nu}_{n}\in\{-h,L+h\}$  are defined similarly.

\item If $Z^{\nu}_{n},Z^{\nu'}_{n}\in\{-h,L+h\}$, then with probability $1$,
\begin{align}\notag
(Z^{\nu}_{n+1},Z^{\nu'}_{n+1})=(Z^{\nu}_{n},Z^{\nu'}_{n})+\left(h(-1)^{\one_{\{Z^{\nu}_{n}=L+h\}}}, h(-1)^{\one_{\{Z^{\nu'}_{n}=L+h\}}}\right)
\end{align}
and
\begin{align}\notag
t_{n+1}=t_n,\qquad
 N^{\nu}_{n+1}= N^{\nu}_{n}+1,\qquad&\text{and}\qquad  N^{\nu'}_{n+1}= N^{\nu'}_{n}+1.
\end{align}

\item For every $n\in\N$, set
\begin{align}\notag
X^{\nu}_n\doteq Z^{\nu}_{M^{\nu}_n}\qquad \text{and}\qquad X^{\nu'}_n\doteq Z^{\nu'}_{M^{\nu'}_n},
\end{align}
where $M^{\nu}_n\doteq\max\{m : N^{\nu}_m\le n\}$ and  $M^{\nu'}_n\doteq\max\{m : N^{\nu'}_m\le n\}$.
\end{itemize}

With the above construction $\{X^{\nu}_n\}$ and $\{X^{\nu'}_n\}$ are controlled Markov chains constructed using the optimal feedback controls $\hat\al^{\nu}$ and
$\hat\al^{\nu'}$ respectively, given on the same probability space. Also relationships \eqref{eq:eq905z}--\eqref{num11} are satisfied by
$\Sigma^\nu$ and $\Sigma^{\nu'}$. The above coupling of the two processes gives the joint evolution  $(X^{\nu}(t),X^{\nu'}(t))_{0\le t\le T}$ as follows.  $X^\nu(0)=X^{\nu'}(0)=x_0$ and for every $t\in\mathbb{T}^h$,
\begin{align}\notag
 &(X^{\nu}(t+\Delta^h),X^{\nu'}(t+\Delta^h))- (X^{\nu}(t),X^{\nu'}(t))\\\notag
 &\qquad\qquad  =
  \begin{cases}
    \left(h\one_{\{X^{\nu}(t)\ne L\}},h\one_{\{X^{\nu'}(t)\ne L\}}\right), & \text{w.p. } (h\min\{b^{\nu}(t),b^{\nu'}(t)\}+\sigma^2)/(2\sigma^2), \\
     \left(h\one_{\{X^{\nu}(t)\ne L\}},-h\one_{\{X^{\nu'}(t)\ne 0\}}\right), & \text{w.p. } h(\Delta b(t))^+/(2\sigma^2), \\
        \left(-h\one_{\{X^{\nu}(t)\ne 0\}},h\one_{\{X^{\nu'}(t)\ne L\}}\right), & \text{w.p. } h(\Delta b(t))^-/(2\sigma^2),\\
        \left(-h\one_{\{X^{\nu}(t)\ne 0\}},-h\one_{\{X^{\nu'}(t)\ne 0\}}\right), & \text{w.p. } (-h\max\{b^{\nu}(t),b^{\nu'}(t)\}+\sigma^2)/(2\sigma^2),
\end{cases}
\end{align}
where
\[\Delta b(t)\doteq b^{\nu}(t)-b^{\nu'}(t),\qquad t\in[0,T].
\]
%
%
We also define the corresponding `unconstrained' increment as
\begin{align}\notag
 &(Z^{\nu}(t+\Delta^h),Z^{\nu'}(t+\Delta^h))-(X^{\nu}(t),X^{\nu'}(t))\\
 &\qquad\qquad  \doteq
  \begin{cases}
    \left(h,h\right), & \text{w.p. } (h\min\{b^{\nu}(t),b^{\nu'}(t)\}+\sigma^2)/(2\sigma^2), \\
     \left(h,-h\right), & \text{w.p. } h(\Delta b(t))^+/(2\sigma^2), \\
        \left(-h,h\right), & \text{w.p. } h(\Delta b(t))^-/(2\sigma^2),\\
        \left(-h,-h\right), & \text{w.p. } (-h\max\{b^{\nu}(t),b^{\nu'}(t)\}+\sigma^2)/(2\sigma^2).
\end{cases}\label{eq:eq1235}
\end{align}
\skp
\noi{\bf Bounding $\boldsymbol{W^2_1(\Phi^h(\nu),\Phi^h(\nu'))}$.}
%
Denote $\Delta X(t)\doteq X^{\nu}(t)-X^{\nu'}(t)$. The processes $\Delta B$, $\Delta R$ and $\Delta Y$ are defined similarly.
Note that
$$W^2_1(\Phi^h(\nu),\Phi^h(\nu'))\le \E[|\Delta X|^2_T].$$
We now estimate $\E[|\Delta X|^2_T]$.
Recall that $\Delta X(0)=0$. From \eqref{num11} and Lemma \ref{lem_Skorohod},
\begin{align}\notag
|\Delta X|_T\le |\Delta X|_T + |\Delta Y|_T + |\Delta R|_T
\le c_S\Big(\int_0^T|\Delta b(s)|ds+\sigma|\Delta B|_T\Big).
\end{align}
Therefore,
\begin{align}\label{ac01}
\E[|\Delta X|^2_T]\le 2c^2_S\E\Big[\Big(\int_0^T|\Delta b(s)|ds\Big)^2+\sigma^2 |\Delta B|^2_T\Big].
\end{align}
We now estimate the second term on the right side. By using the martingale property
of $B^{h,\nu}(t)- B^{h,\nu'}(t)$ and
Doob's inequality,
\begin{align}\label{ac02}
\E[|\Delta B|_T^2]
\le 4\E\Big[\sum_{s\in \mathbb{T}^h}(\Delta B(s+\Delta^h)-\Delta B(s))^2\Big].
\end{align}
From \eqref{eq:eq907} and \eqref{num14},
\begin{align}\notag
\sigma |\Delta B(s+\Delta^h)- \Delta B(s)|&\le\left|(\Delta X+\Delta R-\Delta Y)(s+\Delta^h)- (\Delta X+\Delta R-\Delta Y)(s)\right|+\Delta^h |\Delta b(s)|.
\end{align}
If $(Z^{\nu}(t+ \Delta^h) - X^{\nu}(t))(Z^{\nu'}(t+ \Delta^h) - X^{\nu'}(t)) >0$, i.e.,
 the unconstrained increments are of the same sign, then
\begin{align}\notag
&\left|(\Delta X+\Delta R-\Delta Y)(s+\Delta^h)- (\Delta X+\Delta R-\Delta Y)(s)\right|\\\notag
&\quad=
\left|(X^{\nu}+R^{\nu}-Y^{\nu})(s+\Delta^h)- (X^{\nu}+R^{\nu}-Y^{\nu})(s) \right.\\\notag
&\qquad\qquad\left.-[(X^{\nu'}+R^{\nu'}-Y^{\nu'})(s+\Delta^h)- (X^{\nu'}+R^{\nu'}-Y^{\nu'})(s)]\right|
=0.
\end{align}
If the signs are different, i.e., $(Z^{\nu}(t+ \Delta^h) - X^{\nu}(t))(Z^{\nu'}(t+ \Delta^h) - X^{\nu'}(t)) <0$, then
\begin{align}\notag
&\left|(\Delta X+\Delta R-\Delta Y)(s+\Delta^h)- (\Delta X+\Delta R-\Delta Y)(s)\right|\le 2h.
\end{align}
Hence,
\begin{align}\label{1031b}
\sigma|(\Delta B(s+\Delta^h)-\Delta B(s))|\le 2h\one_{\hat E_s}+|\Delta b(s)|\Delta^h,
\end{align}
where
\begin{align}\notag
\hat E_s\doteq\left\{\omega : (Z^{\nu}(t+ \Delta^h) - X^{\nu}(t))(Z^{\nu'}(t+ \Delta^h) - X^{\nu'}(t)) <0\right\}.
\end{align}
From \eqref{eq:eq1235} we now have that,
\begin{align}\notag
\PP(\hat E_s|\calH^{h}_{s})\le h|\Delta b(s)|/(2\sigma^2),
\end{align}
where $\{\calH^{h}_t\}$ is the filtration generated by the process $(X^{\nu}(t),X^{\nu'}(t))_{0\le t\le T}$.
As a consequence,
\begin{align}\label{1032}
\E\Big[\Big(\sigma(\Delta B(s+\Delta^h)-\Delta B(s))\Big)^2\mid\calH^{h}_{s}\Big]
&\le\E\Big[\left(
2h\one_{\hat E_s}+|\Delta b(s)|\Delta^h
\right)^2\mid\calH^{h}_{s}\Big]\\\notag
&\le2h\,|\Delta b(s)|(\Delta^h+Ch\Delta^h),
\end{align}
where in the above expression, and in the rest of the proof, $C$ refers to a finite positive constant that is independent of $h$ and $s, \nu, \nu'$ and which  can change from one line to the next. Applying the above bound to \eqref{ac02} and taking $h$ sufficiently small such that $Ch\le 1/2$, we get that for sufficiently small $h$, 
\begin{align}\label{1033a}
\sigma^2\E[|\Delta B|^2_T]&\le 12\Delta^h\sum_{s\in\mathbb{T}^h}h\E[|\Delta b(s)|]=12\E\Big[\int_0^Th|\Delta b(s)|ds\Big].
\end{align}
Combining this with \eqref{ac01} and using the inequality,
\begin{align}\label{1025}
h\int_0^T|\Delta b(s)|ds
&\le\frac{1}{2} \Big[h^2T^{1/2}+T^{-1/2}\Big(\int_0^T|\Delta b(s)|ds\Big)^2\Big]
\le \frac{1}{2}T^{1/2}\Big[h^2 + \int_0^T|\Delta b(s)|^2ds\Big]
\end{align}
we get that
\begin{align}\label{ac03}
\E[|\Delta X|^2_T]&\le C\E\Big[\Big(\int_0^T|\Delta b(s)|ds\Big)^2+h\int_0^T |\Delta b(s)|ds\Big]\\\notag
& \le CT^{1/2}h^2+C(T^{1/2}+T)\E\Big[\int_0^T(\Delta W(s))^2ds+\int_0^T(\Delta X(s))^2ds+\int_0^T(\Delta \calD^{h}_xV^{h}(s))^2ds
 \Big],
\end{align}
where for $t\in[0,T]$
\[
\Delta \calD^{h}_xV^{h}(t)\doteq \calD^{h}_xV^{h}_{\nu}(l^h(t),X^\nu(t))-\calD^{h}_xV^{h}_{\nu'}(l^h(t),X^\nu(t)), \, \Delta W(t) \doteq W_1(\nu(l^h(t)),\nu'(l^h(t)))
\]
and
the above inequality also uses the Lipschitz property of $b$ (Assumption \ref{assumption1}), the Lipschitz property of $\hat\al$ (Assumption \ref{assu:lips_al}) and \eqref{1006d}.

We now consider the last term on the right side of \eqref{ac03}. We will
 show that for some $\hat T$ that does not depend on $h, \nu, \nu'$ and all $T \le \hat T$,
\begin{align}\label{1026}
&\E \Big[\int_0^T(\Delta \calD^{h}_xV^h(s))^2ds\Big]\le C\left(h^2+ \sup_{0\le s\le T}(\Delta W(s))^2
+\E[|\Delta X|_T^2]\right).
\end{align}

Define,
\begin{align}\notag
\Delta V^{h}(s)&\doteq V^{h}_\nu(s,X^{\nu}(s))-V^{h}_{\nu'}(s,X^{\nu'}(s)),\;
\Delta g(T)\doteq \Delta V^{h}(T) = g(\nu(T),X^{\nu}(T))-g(\nu'(T),X^{\nu'}(T)),\notag\\
\Delta f(s)&\doteq f(l^h(s),\nu(l^h(s)),X^{\nu}(s),\hat\al^{h,\nu}(s,X^{\nu}(s)))-f(l^h(s),\nu'(s),X^{\nu'}(l^h(s)),\hat\al^{h,\nu'}(s,X^{\nu'}(s))).\notag
\end{align}
From \eqref{1006b} we get that
\begin{align}\notag
&\Delta V^{h}(0)+ \sigma\sum_{s\in\mathbb{T}^h}\Delta \calD^{h}_xV^h(s)(B^{\nu}(s+\Delta^h)-B^{\nu}(s))\\\notag
&\quad=- \sigma\sum_{s\in\mathbb{T}^h}\calD^{h}_xV^{h}_{\nu'}(s,X^{\nu'}(s))(\Delta B(s+\Delta^h)-\Delta B(s))+\Delta g(T)+\int_0^T\Delta f(s)ds\\\notag
&\;\qquad+\int_0^Ty(s)d(\Delta Y(s))+\int_0^Tr(s)d(\Delta R(s)).
\end{align}
By squaring both sides and taking expectations
, we have  
\begin{align}\label{1029}
&(\Delta V^{h}(0))^2+ \E\Big[\Big(\sigma\sum_{s\in\mathbb{T}^h}\Delta\calD^{h}_xV^{h}(s)(B^{\nu}(s+\Delta^h)-B^{\nu}(s))\Big)^2\Big]\\\notag
&\quad\le 5\left\{\E\Big[\Big(\sigma\sum_{s\in\mathbb{T}^h}\calD^{h}_xV^{h}_{\nu'}(s,X^{\nu'}(s))(\Delta B(s+\Delta^h)-\Delta B(s))\Big)^2\Big]+\E[(\Delta g(T))^2]\right.
\\\notag
&\qquad\qquad+\left.\E\Big[\Big(\int_0^T\Delta f(s)ds\Big)^2\Big]+  \E\Big[\Big(\int_0^Ty(s)d(\Delta Y(s))\Big)^2\Big]+\E\Big[\Big(\int_0^Tr(s)d(\Delta R(s))\Big)^2\Big]\right\}.
\end{align}
Here we have used the fact that $B^{h,\nu}$ is a $\{\calG^{h,\nu}_t\}$-martingale  and therefore has mean $0$. Also, we use the elementary inequality that $(\sum_{i=1}^5 a_i)^2 \le 5 \sum_{i=1}^5 a_i^2$.
We now bound the terms in the inequality above. For any $s\in\mathbb{T}^h$, we get by \eqref{num5} and \eqref{eq:eq907} that
\begin{align}\notag
\E\Big[\Big(\Delta\calD^{h}_xV^{h}(s)(B^{\nu}(s+\Delta^h)-B^{\nu}(s)\Big)^2\mid\calH^{h}_{s}\Big]=(\Delta\calD^{h}_xV^{h}(s))^2(\Delta^h+C(\Delta^h)^2).
\end{align}
Using once more the martingale property of $B^{\nu}$
 we get that, for  sufficiently small $h$,
\begin{align}
\E\Big[\Big(\sum_{s\in\mathbb{T}^h}\Delta\calD^{h}_xV^{h}(s)(B^{\nu}(s+\Delta^h)-B^{\nu}(s))\Big)^2\Big]
&=\sum_{s\in\mathbb{T}^h}\E\Big[\Big(\Delta\calD^{h}_xV^{h}(s)(B^{\nu}(s+\Delta^h)-B^{\nu}(s))\Big)^2\Big]\notag\\
&\ge\frac{1}{2}\E\Big[\int_0^T(\Delta \calD^{h}_xV^{h}(s))^2ds\Big].\label{1031}
\end{align}

Recall that from \eqref{1006c}, $|\calD^{h}_xV^{h}_{\nu'}(\cdot,X^{\nu'}(\cdot))|_T\le c_d(M)$ whenever $T \le M$. Henceforth we will only consider $T \in [0,M]$.
Using \eqref{1032}, we get that for sufficiently small $h$,
\begin{align}\label{1033}
&\E\Big[\Big(\sigma\sum_{s\in\mathbb{T}^h}\calD^{h}_xV^{h}_{\nu'}(s,X^{\nu'}(s))(\Delta B(s+\Delta^h)-\Delta B(s))\Big)^2\Big]\\\notag
&\quad=\sum_{s\in\mathbb{T}^h}\E\Big[\Big(\sigma\calD^{h}_xV^{h}_{\nu'}(s,X^{\nu'}(s))(\Delta B(s+\Delta^h)-\Delta B(s))\Big)^2\Big]\\\notag
&\quad\le 3(c_d(M))^2\E\Big[\int_0^Th|\Delta b(s)|ds\Big]\\\notag
&\quad \le CT^{1/2}\E\Big[h^2+\int_0^T(\Delta W(s))^2ds+\int_0^T(\Delta X(s))^2ds+\int_0^T(\Delta \calD^{h}_xV^{h}(s)))^2ds
 \Big],
\end{align}
where the last inequality follows by a similar inequality as in \eqref{1025}.

From \eqref{19b}, it is easy to see that
\begin{align}\label{1034}
\E[(\Delta g(T))^2]&\le C\left((\Delta W(T))^2+\E[(\Delta X (T))^2]\right),
\end{align}
and that
\begin{align}
\E\Big[\left(\int_0^T\Delta f(s)ds\right)^2\Big]
&\le
C\E\Big\{\Big(\int_0^T\Delta W(s)ds\Big)^2+ \Big(\int_0^T\Delta X(s)ds\Big)^2 +\Big(\int_0^T\Delta\calD^{h}_xV^{h}(s)ds\Big)^2\Big\}\notag\\
&\le
CT\E\Big[\int_0^T(\Delta W(s))^2ds+\int_0^T(\Delta X(s))^2ds+\int_0^T(\Delta \calD^{h}_xV^{h}(s)))^2ds
 \Big].\label{1035}
\end{align}
Using integration by parts and the boundedness of $y$
%
\begin{align}\label{1037}
\E\Big[&\Big(\int_0^Ty(s)d(\Delta Y(s))\Big)^2\Big]\\\notag
&\le C\E\Big[ |\Delta Y|^2_T
\Big]\\\notag
&\le C\E\Big[|F^{\nu}+\sigma B^{\nu}-F^{\nu'}-\sigma B^{\nu'}|^2_T\Big]\\\notag
&\le C\E\Big[\Big(\int_0^T|\Delta b(s)|ds\Big)^2+\sigma^2|\Delta B|_T^2\Big]\\\notag
& \le CT^{1/2}h^2+C(T^{1/2}+T)\E\Big[\int_0^T(\Delta W(s))^2ds+\int_0^T(\Delta X(s))^2ds+\int_0^T(\Delta \calD^{h}_xV^{h}(s)))^2ds
 \Big],
\end{align}
where the inequality on the third line is from Lemma \ref{lem_Skorohod}, the fourth line is from \eqref{num14} and the last line uses \eqref{1033a} and \eqref{ac03}.
A similar bound holds for $\E\Big[\Big(\int_0^Tr(s)d(\Delta R(s))\Big)^2\Big]$.

From \eqref{1029}--\eqref{1037}, we have
\begin{align}\label{1038}
&(\Delta V^{h}(0))^2+ \frac{1}{2}\sigma^2\E\Big[\int_0^T(\Delta \calD^{h}_xV^{h}(s))^2ds\Big]\\\notag
&\quad\le C\left( (\Delta W(T))^2+\E[(\Delta X(T))^2]\right) \\\notag
&\qquad+CT^{1/2}h^2+C(T^{1/2}+T)\E\Big[\int_0^T(\Delta W(s))^2ds+\int_0^T(\Delta X(s))^2ds+\int_0^T(\Delta \calD^{h}_xV^{h}(s)))^2ds
 \Big].
\end{align}
Thus we can find a $\hat T_1 \in (0,M)$ and $\hat h_1>0$ such that for all $T \le \hat T_1$ and $h \le \hat h_1 \wedge \hat T_1$,
 \eqref{1026} is satisfied. Together with \eqref{ac03} we now get that there exist $\hat T  \in (0,M)$, $\hat h>0$ and $q\in(0,1)$ such that for every $T\le \hat T$ and   $h\in(0,\hat h \wedge \hat T)$,
\begin{align}\notag
W^2_1(\Phi^h(\nu),\Phi^h(\nu'))
&\le\E[|\Delta X|^2_T]\le q\Big(h^2+\sup_{0\le s\le T}(\Delta W(s))^2\Big)\le q\Big(h^2+W^2_1(\nu,\nu')\Big).
\end{align}


\hfill$\Box$
\skp

We finally prove the last statement in Lemma \ref{lem_52}, namely the inequality in \eqref{1006c}.\\

\noi {\bf Proof of \eqref{1006c}.} Fix $\nu\in\calPTL$, $x\ne x'$ in $\mathbb{S}^h_0$, and $t_0\in\mathbb{T}^h$, which will be regarded as the initial time. As in the proof of Proposition \ref{prop_h_cont}, one can define a coupling of two processes on the same $h$-grid, both of which are driven by the same $\nu$. The first one is denoted as $(X(s),Y(s),R(s),B(s),\al(s))_{t_0\le s\le T}$, where its components are defined in \eqref{eq:eq905z}--\eqref{eq:eq907} with $X(t_0)=x$ and $\al$ is the optimal policy for this process. The second process, denoted as $(X'(s),Y'(s),R'(s),B'(s),\al(s))_{t_0\le s\le T}$ is also given by \eqref{eq:eq905z}--\eqref{eq:eq907} using the same control process $\{\alpha(s)\}$ as for the first one, except that the second one  starts at $x'$, i.e., $X'(t_0)=x'$.

For every $s\in[t_0,T]$, let $\Delta X(s)\doteq X(s)-X'(s)$ and
\begin{align}\notag
\Delta b(s)\doteq b(l^h(s),\nu(l^h(s)),X(s),\al(s))-b(l^h(s),\nu(l^h(s)),X'(s),\al(s)).
\end{align}
Processes $\Delta Y(s)$ and $\Delta R(s)$ are defined in a similar manner.
By definition $\Delta X(0)=x-x'$. The arguments that lead to \eqref{ac03} can also be applied here, but in fact they are simpler here since $\Delta W=0$ and the controls are the same for both processes.
Specifically,
\begin{align}\label{ac05}
&\E\Big[ \sup_{t_0\le s\le T} \left(|\Delta X(s)| + |\Delta R(s)| + |\Delta Y(s)|\right)^2\Big]\\\notag
&\quad\le C(x-x')^2+C_1(T-t_0)^{1/2}h^2+C((T-t_0)^{1/2}+(T-t_0))\E\Big[\int_{t_0}^T(\Delta X(s))^2ds
 \Big]\\\notag
&\quad\le C(x-x')^2(1+T^{1/2})+C(T^{1/2}+T)\E\Big[\int_{t_0}^T(\Delta X(s))^2ds
 \Big].
\end{align}
 The second inequality is a consequence of the fact that since $x\ne x'$, $h\le |x-x'|$.
Therefore,
\begin{align}\notag
&\E\Big[ \sup_{t_0\le s\le T}(\Delta X(s))^2\Big]
\le C(x-x')^2(1+T^{1/2})+C(T^{1/2}+T)\E\Big[\int_{t_0}^T(\Delta X(s))^2ds
 \Big].
\end{align}
By 
Gr{\"o}nwall's inequality,
$\E\Big[ \sup_{t_0\le s\le T}(\Delta X(s))^2\Big] \le c_T|x-x'|^2$,
where $$c_{T}\doteq C(1+T^{1/2})\exp\left\{C_1(T^{3/2}+T^2)\right\}.$$ Using the above bound  we have that the left side in \eqref{ac05} can be bounded above by
\begin{align}\notag
 C(x-x')^2(1+T^{1/2})+C(T^{1/2}+T)T\;
\E\Big[\sup_{t_0\le s\le T}(\Delta X(s))^2
 \Big]
 \le \tilde c_{T} |x-x'|^2,
\end{align}
where $\tilde c_{T}:=\left(4CL^2(1+T^{1/2})+C(T^{1/2}+T)Tc_{T}\right)^{1/2}$, and so,
\begin{align}\notag
\E\Big[ \sup_{t_0\le s\le T} \left( |\Delta X(s)| + |\Delta Y(s)| + |\Delta R(s)| \right)\Big]\le (\tilde c_{T})^{1/2} |x-x'|,
\end{align}
Consequently, using integration by parts as in \eqref{1037} and Lipschitz property of $f$ and $g$, we get that
\begin{align}\notag
&V^{h}_\nu(t_0,x')-V^{h}_\nu(t_0,x)\\\notag
&\quad\le J^{h,\nu}(t_0,x',\al)-J^{h,\nu}(t_0,x,\al)\\\notag
&\quad\le\E\Big[|g(\nu(T),X'(T))-g(\nu(T),X(T))|+ \left|\int_{t_0}^Ty(s)d(\Delta Y(s))\right|+
\left|\int_{t_0}^Tr(s)d(\Delta Y(s))\right|
\\\notag
&\qquad\qquad +\int_{t_0}^T\left|f(l^h(s),\nu(l^h(s)),X'(s),\al(s))-f(l^h(s),\nu(l^h(s)),X(s),\al(s))\right|ds\Big]\\\notag
&\quad\le \bar c_T|x-x'|,
\end{align}
where $\bar c_T$ depends on the parameters $c_L,c_S$, the bounds on $y$ and $r$,  and the terminal time $T$.
By reversing the roles of the processes we get that for every $t_0\in\mathbb{T}^h$, 
\begin{align}\notag
|V^{h}_\nu(t_0,x)-V^{h}_\nu(t_0,x')|\le \bar c_T|x-x'|
\end{align}
and the result follows.
\hfill$\Box$

\skp\noi

 \section{Numerical study} \label{sec:numerics}
In this section we present a numerical example. 
We set the parameters $L=1$, $T=0.4$, $\sigma=1$, and $U=\{-0.75,0.25\}$. Also, $b(t,\eta,x,\al)=2x+7\al$, $f(t,\eta,x,\al)=(4x-5\bar\eta)^2+\al^2$, $g(\eta,x)=(4x-5\bar\eta)^2$, $y(t,x)=0$, and $r(t,x)=15$, where $\bar\eta$ is the mean of $\eta$. Assumption \ref{assumptionU} obviously holds and therefore, by Proposition \ref{thm_fixed} the MFG admits a unique solution and Assumption  \ref{assu:uniqmfg} holds. 
 The initial state of the MFG is taken to be $X(0)=0.5$ and in the numerical scheme $x^{h}=X^{h}(0)=\lfloor x/h \rfloor\cdot h$. We choose the initial function $\nu^{h}(t)=\delta_{\lfloor x/h \rfloor\cdot h}$, $t\in[0,T]$.

\begin{figure}[h!]
 \centering
 \includegraphics[width=1\textwidth]{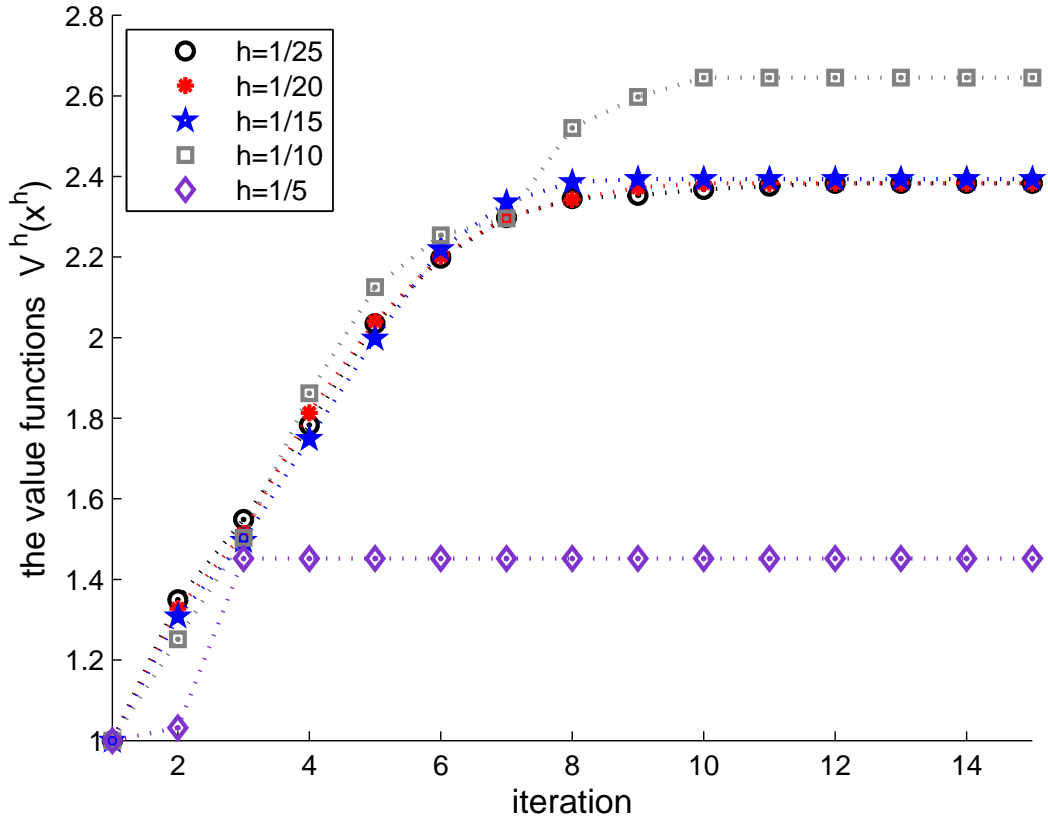}
 \caption{The value functions $V^h(x^h)$ for the numerical example.}
 \label{figure1}
 \end{figure}

\begin{figure}[h!]
 \centering
 \includegraphics[width=1\textwidth]{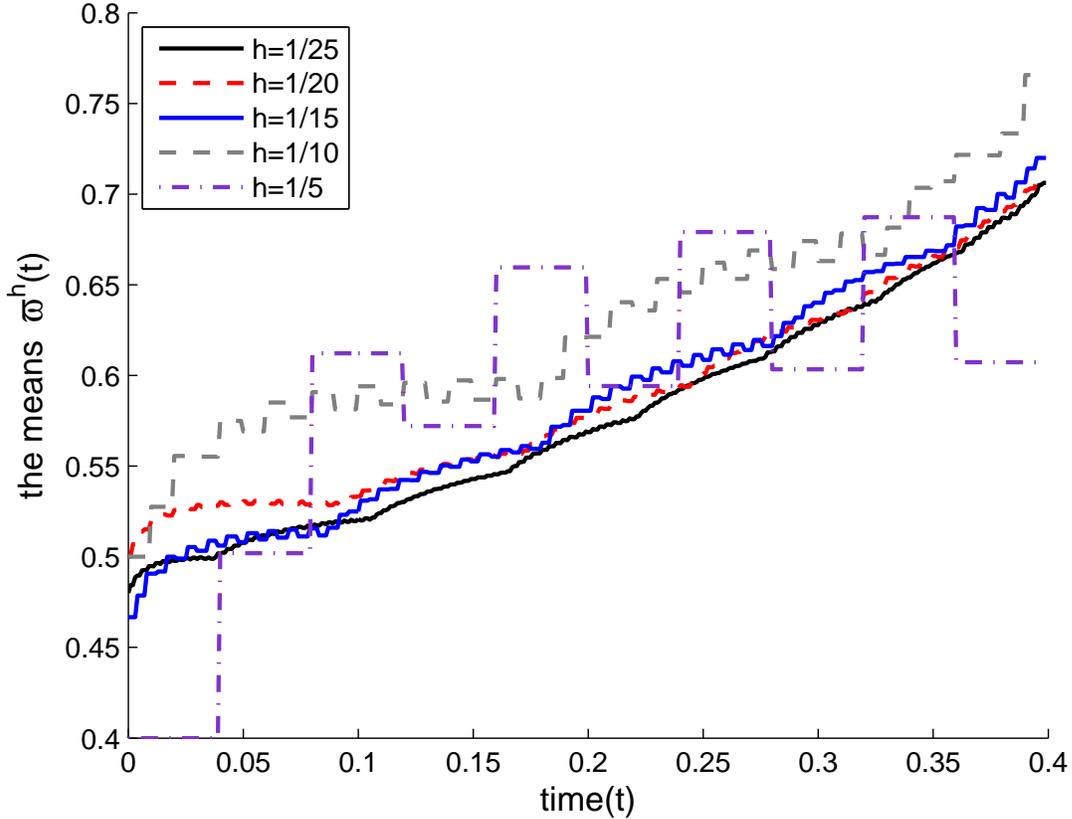}
 \caption{The means $\varpi^{h}$ for the numerical example.}
 \label{figure2}
 \end{figure}

\begin{figure}[h!]
 \centering
 \includegraphics[width=1\textwidth]{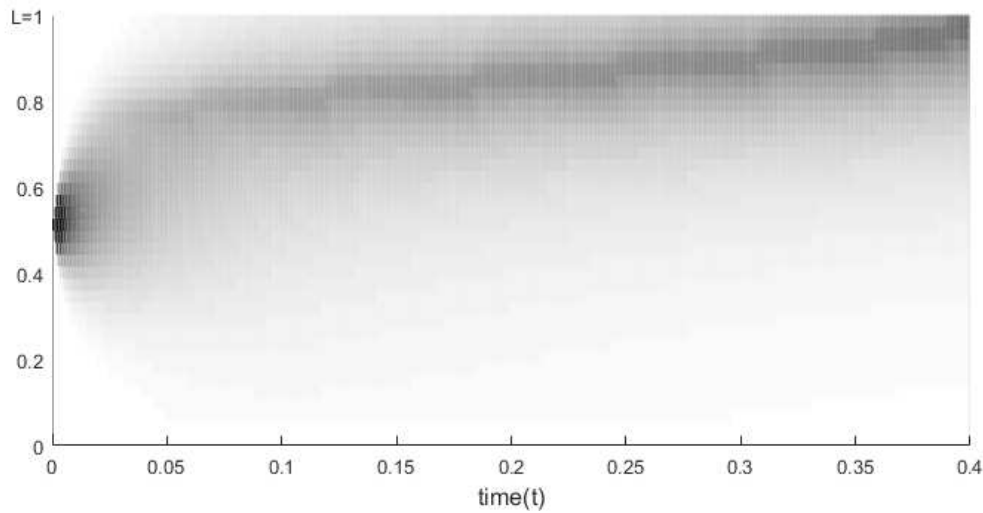}
 \caption{Bird's eye view of $\nu^{h}$ for the numerical example with $h=25$.}
 \label{figure3}
 \end{figure}

 We implemented the algorithm described in Construction \ref{const:const2} by computing $15$ iterations of the map $\Phi^h$ for each of the $h$'s taken from the set $\{1/5,1/10,1/15,1/20,1/25\}$. For each $h$ we calculated the value function of the MDP after each iteration $V^{h}(x^{h})$. 
 Since our example depends on $\nu^{h}$ only through its mean, we also calculated the mean of $\nu^{h}(\cdot)$, which we denote by $\varpi^{h}$.
Figure \ref{figure1} illustrates the convergence of the value functions of the MDP's to the value function of the MFG. The convergence of the means $\varpi^{h}$ is illustrated in Figure \ref{figure2}. Finally in Figure \ref{figure3} we present the distribution taken from the last iteration, $\nu^{1/25}$. Here we provide a bird's eye view of this distribution, where the darker areas represent greater density.

 \renewcommand{\arraystretch}{1.3}

%


 As can be seen by Figure \ref{figure2}, shortly after time zero, the means tend to increase w.r.t.~the time. The reason for this is that there are two opposing ``forces'' in our example. The reflection cost $r=15$ ``pushes'' the process $X^{h}$ downwards when it is close to the boundary $L$ and then the optimal control  $\al$ is close to $-0.75$. When $X^{h}$ is relatively far away from the boundary, the control $\al=-0.75$ is too costly and then the optimal control $\al$ is close to $0.25$. As we approach the terminal time, the rejection cost has less impact and therefore the distribution has higher expectation. 

\skp\noi{\bf Acknowledgement.} We are thankful to the anonymous referee and AE for their suggestions, which helped us to improve the presentation of the paper.

\footnotesize{
\bibliographystyle{abbrv} 
\bibliography{bib_Asaf_IMS} 
}\footnotesize{{\sc
\bigskip

\noindent
Erhan Bayraktar\\
Department of Mathematics\\
University of Michigan\\
Ann Arbor, MI 48109, USA\\
email: erhan@umich.edu\\
web: www-personal.umich.edu/$\sim$ erhan/
\skp

\noindent
Amarjit Budhiraja \\
Department of Statistics and Operations Research\\
University of North Carolina\\
Chapel Hill, NC 27599, USA\\
email: budhiraj@email.unc.edu
web: http://www.unc.edu/$\sim$ budhiraj/
\skp

\noindent
Asaf Cohen\\
Department of Statistics\\
University of Haifa\\
Haifa 31905, Israel\\
email: shloshim@gmail.com\\
web: https://sites.google.com/site/asafcohentau/

}}
\end{document}